\documentclass[a4paper,11pt,oneside]{article}

\usepackage{abstract}
\usepackage{amsmath,amsthm,amsfonts,amssymb,amscd}
\usepackage{graphicx}
\usepackage{mathtools}
\usepackage{hyperref}
\usepackage{cleveref}
\usepackage{mathrsfs}
\usepackage{float}
\usepackage{csvsimple}
\usepackage{enumerate}
\usepackage{booktabs}
\newtheorem{remark}{Remark}

\newtheorem{Algorithm}{Algorithm}
\newtheorem{problem}{Problem}
\newtheorem{theorem}{Theorem}
\usepackage{pgfplots}
\pgfplotsset{compat=1.14}
\usepackage[noend]{algpseudocode}
\usepackage{graphicx,bm,xcolor}
\usepackage{algorithm}
\usepackage{pdfpages}
\usepackage{algpseudocode}
\usepackage{stmaryrd}
\usepackage{hyperref}
\usetikzlibrary{arrows}
\usepackage{caption}
\usepackage[T1]{fontenc}
\usepackage{t1enc}
\usepackage[ngerman,english]{babel}
\usepackage{verbatim} 

\usepackage{url}

\usepackage{units}  
\usepackage{siunitx}

\usepackage{tabularx}
\usepackage[font={footnotesize},%
  labelsep=colon]{caption}
\usepackage{tabu}
\usepackage{tikz}
\tikzset{every path/.append style={line width=1pt}}

\usepackage{authblk}

\usepackage[top=2cm, bottom=2cm, left=2cm, right=2cm]{geometry}

\begin{document}

\title{\textbf{Mesh adaptivity for quasi-static phase-field fractures based on a residual-type a posteriori error estimator}}

\author[1]{Katrin Mang}
\author[2]{Mirjam Walloth}
\author[1]{Thomas Wick}
\author[2]{Winnifried Wollner}

\affil[1]{Leibniz Universit\"at Hannover, Institut f\"ur Angewandte
  Mathematik\\
  
AG Wissenschaftliches Rechnen, Welfengarten 1, 30167 Hannover, Germany}
\affil[2]{Department of Mathematics, Technische Universit{\"a}t Darmstadt\\

Dolivostrasse 15, 64293 Darmstadt, Germany}

\date{\textit{\today}}

\maketitle

\renewenvironment{abstract}
{\begin{quote}
\noindent \rule{\linewidth}{.5pt}\par{\bfseries \abstractname.}}
{\medskip\noindent \rule{\linewidth}{.5pt}
\end{quote}
}

\begin{center}
 \textbf{This is the preprint version  of an accepted article to be published in
the GAMM-Mitteilungen 2019\\ \url{https://onlinelibrary.wiley.com/journal/15222608}}
\end{center}

\begin{abstract}
In this work, we consider adaptive mesh refinement for
a monolithic phase-field description for fractures
in brittle materials. Our approach is based on an a posteriori error estimator
for the phase-field variational 
inequality realizing the fracture irreversibility constraint. 
The key goal is the development of a reliable and efficient residual-type error estimator 
for the phase-field fracture model in each time-step.
Based on this error estimator, 
error indicators for local mesh adaptivity are extracted.
The proposed estimator is based on a technique known for 
singularly perturbed equations
in combination with estimators for variational inequalities.
These theoretical developments are used to formulate an adaptive mesh 
refinement algorithm. For the numerical solution, 
the fracture irreversibility is imposed using a Lagrange multiplier. The resulting 
saddle-point system has three unknowns: displacements, phase-field, 
and a Lagrange multiplier for the crack irreversibility.
Several numerical experiments demonstrate our theoretical findings
with the newly developed estimators and the corresponding refinement strategy.
\end{abstract}

\section{Introduction}
\label{sec_intro}
Fracture propagation and damage mechanics are current topics in 
theoretical mathematics (calculus of variations), numerical 
mathematics, and engineering. A well-established variational approach for Griffith's~\cite{griffith1920phenomena} 
quasi-static brittle fracture was introduced by Francfort
and Marigo~\cite{FraMar98}. 
The focus of the current work is on the development of 
a posteriori error estimation 
and local mesh adaptivity for such a variational (phase-field) fracture 
formulation. 

Important motivations for mesh adaptivity are, first, the 
challenge of the 
resolution of the phase-field regularization parameter 
$\epsilon$ in relation 
to the spatial discretization parameter $h$ such that $h\ll\epsilon$. 
Second, in many applications, the crack tip is of specific interest. Thus 
having an error estimator along with localized crack tip mesh refinement 
would be beneficial. We notice that mesh refinement for phase-field fracture 
problems may be problematic since it may occur that the (unknown) fracture path 
depends on the locally refined mesh. This yields numerical solutions purely 
dependent on the mesh refinement algorithm (see e.g.,~\cite{ArFoMiPe15}) 
and discussions and references provided therein.

The first study on local mesh adaptivity for phase-field fracture
was undertaken in~\cite{BuOrSue10}. An extension to
anisotropic mesh adaptivity was done in~\cite{ArFoMiPe15}.
Goal-oriented error estimation using 
dual-weighted residuals was addressed in~\cite{Wi16_dwr_pff}.
Another method (only mesh refinement, but no error estimator) 
that purely focuses on fine meshes 
in the crack region has been developed in~\cite{HeWheWi15}
for simulations in two spatial dimensions 
and the extension for three-dimensional computations was 
considered in~\cite{LeeWheWi16}. Using these last developments, 
a computational convergence analysis using high performance parallel computing 
and local mesh adaptivity was carried out in~\cite{HeiWi18_pamm}.
All these studies show that local mesh refinement is a key ingredient 
for phase-field fractures, in particular in view of 
working with sufficiently small phase-field regularization parameters.
However, due to the complexity of the problem setting such as 
nonlinearities and variational inequalities (due to the crack 
irreversibility constraint), further studies are necessary.

It has been observed, e.g., in~\cite{BuOrSue10}, that
  asserting convergence in the individual residuals in each time step
  is sufficient to
  obtain stationary points of the phase-field problem. These residuals
correspond to linear elasticity with a degenerated coefficient, and a
singularly perturbed obstacle problem. 
Examples of linear elliptic problems with inequality constraints 
are obstacle and contact problems~\cite{Rodrigues:1987,KikuchiOden:1988}. 
Different types of error estimators for obstacle problems can be found in, e.g.,~\cite{ChenNochetto:2000, WeissWohlmuth:2010, ZouVeeserKornhuberGraeser:2011}. 
By measuring the error in the solutions as well as in the constraining
force, 
the first efficient and reliable residual-type estimator for obstacle problems 
has been derived in~\cite{Veeser:2001}. This approach has been
extended 
to discontinuous Galerkin methods in~\cite{GudiPorwal:2014, GudiPorwal:2016}. 
Error estimators for contact problems are given 
in~\cite{HildNicaise:2007, WeissWohlmuth:2009, KrauseVeeserWalloth:2015}. 
In~\cite{MoonNochettoPetersdorffZhang:2007, KrauseVeeserWalloth:2015} 
the local structure of the solution and constraining force has been
exploited 
to localize the estimator contribution related to the constraints. This approach
enables a good resolution of the critical regions between the phases of
active 
and non-active constraints while avoiding over-refinement in the active set. 
An extension of this approach to discontinuous Galerkin methods has
been used in~\cite{Walloth:2019}.

The main contributions of this work are concerned with 
the development of an adaptive refinement strategy 
based on an a posteriori error estimator for the phase-field variational
inequality, with particular emphasis on the robustness of the
estimator with respect to the phase-field parameter $\epsilon$.
This aspect is a first difference to previous studies, 
in particular in comparison to~\cite{ArFoMiPe15} and~\cite{BuOrSue10}
in which the irreversibility was modeled by enforcing zero values along
the fracture and robustness of the estimator was not considered.
A second difference to previous works is that our adaptive refinement
  procedure will be based on consecutive solutions over the whole time
  interval on relatively coarse meshes, resulting in refinement indicators
for the entire time horizon.
These theoretical and algorithmic derivations are substantiated 
with appropriate numerical tests including studies with varying the relationship of the
discretization parameter $h$ and the crack width $\epsilon$. In
  one numerical example, we are concerned with very fast crack growth in order to 
  study whether the proposed estimator can cope with such situations.

The outline of this paper is as follows: In Section~\ref{sec_modeling},
we introduce the notation and our model formulation. Next, in Section~\ref{sec_discretization}, the discretization and the a posteriori error
estimator for the phase-field variational inequality
are developed. We complement this by providing the details for our
  numerical simulation in Section~\ref{sec:solver}.
In Section~\ref{sec_tests}, some numerical tests are carried
out for showing 
the performance of our theoretical developments (more details are provided in the last chapters at the end). We summarize our 
observations in Section~\ref{sec_conclusions}.


\section{Modeling}\label{sec_modeling}
This section covers the basic notation, the function spaces, the
required variables and the strong and weak problem formulation of a
phase-field approach. We emanate from a two-dimensional, polygonal
domain $\Omega\subset \mathbb{R}^2$. By means of an elliptic
functional developed by Ambrosio-Tortorelli~\cite{ambrosio1992approximation,ambrosio1990approximation} the
lower-dimensional crack $C \subset \Omega$ is approximated by a
phase-field variable $\varphi: (\Omega \times I) \to [0,1]$ with
$\varphi=0$ in the crack and $\varphi=1$ in the unbroken material.
A parameter $\epsilon>0$ determines the width of a transition zone
between the unbroken material and the broken material inside the
approximate crack.

Let $I$ be a loading (time) interval $[0,T]$, where $T>0$ is the end time value.
A displacement function $\boldsymbol{u}:(\Omega \times I) \to \mathbb{R}^2$ is defined on the domain $\Omega$.
The boundary $\Gamma=\partial\Omega$ is a Dirichlet boundary for the displacements $\boldsymbol{u}$. 
For the phase-field variable, we have Neumann values $\boldsymbol{\nabla}\varphi\cdot \boldsymbol{n}= 0$ on the whole boundary $\Gamma$ where $\boldsymbol{n}$ is the unit outward normal to the boundary.
The physics of the underlying problem ask to enforce that the fracture
cannot heal, i.e., that $\varphi$ is monotone non-increasing with respect
to $t\in I$. This condition is called irreversibility condition.

\subsection{Strong formulation} 

In order to give the strong formulation of our model problem, we need some further definitions. 
The Frobenius scalar product of two matrices of the same dimension is denoted as $(\boldsymbol{A}:\boldsymbol{B}):= \sum_i \sum_j a_{ij} b_{ij}$. 

A degradation function $g(\varphi)$ is defined as
\begin{align*}
g(\varphi):=(1-\kappa)\varphi^2 + \kappa,
\end{align*}
with a small regularization parameter $\kappa > 0$.
The stress tensor $\boldsymbol{\sigma}(\boldsymbol{u})$ is given by
\begin{align*}
\boldsymbol{\sigma}(\boldsymbol{u}) := 2 \mu \boldsymbol{E}_{\text{lin}}(\boldsymbol{u}) + \lambda \operatorname{tr} (\boldsymbol{E}_{\text{lin}}(\boldsymbol{u})) \textbf{I},
\end{align*}
with the Lam\'e parameters $\mu,\lambda > 0$.
Here, $\boldsymbol{E}_{\text{lin}}(\boldsymbol{u})$ is the linearized strain tensor:
\begin{align*}
\boldsymbol{E}_{\text{lin}}(\boldsymbol{u}):=\frac{1}{2} (\nabla \boldsymbol{u} + \nabla \boldsymbol{u}^T),
\end{align*}
and $\textbf{I}$ denotes the two-dimensional identity matrix.
In Miehe et al.~\cite{MieWelHof10a} a stress splitting has been proposed for fracture phase-field models. 
The linearized strain tensor is decomposed into its tensile and compressive parts, i.e., $\boldsymbol{E}_{\mathrm{lin}} := \boldsymbol{E}^+_{\mathrm{lin}}  + \boldsymbol{E}^-_{\mathrm{lin}}$ with
\begin{align*}
\boldsymbol{E}^+_{\mathrm{lin}} := \boldsymbol{Q}\boldsymbol{D}^+ \boldsymbol{Q}^T,
\end{align*} 
where $\boldsymbol{Q}$ is the matrix of eigenvectors and $\boldsymbol{D}$ the matrix with the eigenvalues on the diagonal. Further, $(\cdot)^+$ denotes the positive part, i.e., on the diagonal of $\boldsymbol{D}^+$ are either the positive eigenvalues or zeros. 
We use the stress splitting of~\cite{MieWelHof10a} which is given by
\begin{align*}
\boldsymbol{\sigma}^+(\boldsymbol{u}) :=& 2\mu\; \boldsymbol{E}^+_{\mathrm{lin}}(\boldsymbol{u}) + \lambda\;\mathrm{max}\{0,\operatorname{tr}(\boldsymbol{E}_{\mathrm{lin}}(\boldsymbol{u}))\}\textbf{I}, \\[3pt]
\boldsymbol{\sigma}^-(\boldsymbol{u}) := &2\mu\; (\boldsymbol{E}_{\mathrm{lin}}(\boldsymbol{u})-\boldsymbol{E}^+_{\mathrm{lin}}(\boldsymbol{u})) + \lambda\;(\operatorname{tr}(\boldsymbol{E}_{\mathrm{lin}}(\boldsymbol{u}))-\mathrm{max}\{0,\operatorname{tr}(\boldsymbol{E}_{\mathrm{lin}}(\boldsymbol{u}))\}\textbf{I}.
\end{align*}
The continuous formulation attributed to Miehe et al.~\cite{MieWelHof10a,miehe2010phase} is given in the following. 
Find $\boldsymbol{u}:(\Omega\times I) \to \mathbb{R}^2$ and $\varphi:(\Omega \times I) \to \mathbb{R}$ such that
\begin{align*}
-\nabla \cdot (g(\varphi) \boldsymbol{\sigma}^+(\boldsymbol{u}))  + \boldsymbol{\sigma}^-(\boldsymbol{u})) = 0\quad \text{in}\ (\Omega \times I),
\end{align*}
\begin{align}
\begin{aligned}
(1- \kappa){\varphi}\boldsymbol{\sigma}^+(\boldsymbol{u}): \boldsymbol{E}_{\text{lin}}(\boldsymbol{u})- \frac{G_c}{\epsilon}(1-\varphi)+ \epsilon G_c \Delta \varphi \geq 0\quad \text{in}\ (\Omega \times I).\label{phi}
\end{aligned}
\end{align}
Herein, $G_c$  is the critical energy release rate.
The crack irreversibility condition is determined by
\begin{align}
 \partial_t \varphi \leq 0\quad \text{in}\ (\Omega \times I),\label{irre}
\end{align}
which has to be cautiously treated in the numerical solution process. 

The system is completed by continuous Dirichlet boundary conditions
\begin{align*}
  \boldsymbol{u}= \boldsymbol{u}_D\quad \text{on}\ (\Gamma\times I),
\end{align*}
for the displacement function $\boldsymbol{u}$, a complementarity
relation between the phase-field equation in~\eqref{phi} and the crack
irreversibility constraint in~\eqref{irre} as
\begin{align*}
\left((1-\kappa)\varphi \boldsymbol{\sigma}^+(\boldsymbol{u}):\boldsymbol{E}_{\mathrm{lin}}(\boldsymbol{u}) - \frac{G_c}{\epsilon}(1-\varphi) + \epsilon G_c\Delta\varphi\right)\cdot\left(\partial_t\varphi\right) &= 0,
\end{align*}
and an initial condition
\begin{align*}
\varphi(x,0)= \varphi_0\quad \text{in}\ (\Omega \times \{0\}).
\end{align*}

\subsection{Time-discrete weak formulation}
From now on, we consider a time discrete formulation on
a fixed subdivision $0 =t_0 < t_1 < \ldots < t_N = T$ of the interval
$I$. We define approximations $(\boldsymbol{u}^n,\varphi^n) \approx
(\boldsymbol{u}(t_n),\varphi(t_n))$ and hence 
the irreversibility condition is given by
$\varphi^{n}\leq\varphi^{n-1}$ for all $n = 1, \ldots, N$.

For the formulation of the time step problems, we introduce the space 
$\boldsymbol{\mathcal{H}}:=\boldsymbol{H}^1(\Omega;\mathbb R^2)$ and let
$\boldsymbol{u}_D^n = \boldsymbol{u}_D(t^n) \in \boldsymbol{\mathcal{H}}\cap \boldsymbol{\mathcal{C}}^0(\Gamma_D)$ be a
continuation of the Dirichlet-data. For the displacement, we define
$\boldsymbol{\mathcal{H}}_{0}:=\{\boldsymbol{w}\in \boldsymbol{H}^1(\Omega)\mid \boldsymbol{w}=\boldsymbol{0}\mbox{ a.e. on }\Gamma_D\}$, 
where $\boldsymbol{w}=\boldsymbol{0}$ on $\Gamma_D$ has to be understood in the usual trace sense.
The phase-field space is $\mathcal{H} = H^1(\Omega)$ with the feasible set
$\mathcal{K}(\varphi^{n-1}):=\{\psi\in \mathcal{H}\mid \psi\leq \varphi^{n-1}\leq 1\}$. 

In the following, we denote the $L^2$-scalar product by
$\left<\cdot,\cdot\right>$ and the duality pairing of $\mathcal{H}$ with its dual by $\left<\cdot,\cdot\right>_{-1,1}$.
For tensor-valued functions $\boldsymbol{A}$ and $\boldsymbol{B}$, of the same dimensions, it
holds $\left<\boldsymbol{A},\boldsymbol{B}\right> := \int_\Omega
\boldsymbol{A} : \boldsymbol{B}\ \mathrm{d}x$.
The $L^2$-norm is denoted by $\|\cdot\|$ and the $H^1$-norm by $\|\cdot\|_1$. 

Now, we can state the problem for calculating the solution $(\boldsymbol{u}^n,\varphi^n)$: 
\begin{problem}[Weak formulation in each time step]\label{WeakFormulation}
For $n=1,\ldots,N$, and given $\varphi^{n-1}$, find $(\boldsymbol{u}^n, \varphi^n)\in(\boldsymbol{\mathcal{H}}_0+\{\boldsymbol{u}^n_D\})\times\mathcal{K}(\varphi^{n-1})$ such that
\begin{equation}\label{VI_orig}
  \begin{aligned}
    \left<g(\varphi^{n-1})\boldsymbol{\sigma}^+(\boldsymbol{u}^n), \boldsymbol{E}_{\mathrm{lin}}(\boldsymbol{w})\right> + \left<\boldsymbol{\sigma}^-(\boldsymbol{u}^n) ,\boldsymbol{E}_{\mathrm{lin}}(\boldsymbol{w}) \right>& = 0\quad \forall \boldsymbol{w}\in \boldsymbol{\mathcal{H}}_0,\\
    \left<(1-\kappa)\varphi^n \boldsymbol{\sigma}^+(\boldsymbol{u}^n):\boldsymbol{E}_{\mathrm{lin}}(\boldsymbol{u}^n),\psi -\varphi^n\right>  - \frac{G_c}{\epsilon}\left<1-\varphi^n, \psi -\varphi^n \right>&\\
    + \epsilon G_c\left<\nabla\varphi^n,\nabla(\psi -\varphi^n)\right>  &\ge 0\quad \forall \psi\in\mathcal{K}(\varphi^{n-1}).
  \end{aligned}
\end{equation}
\end{problem}

\begin{remark}
 In the first term we time-lag the phase-field variable in order to convexify the problem; see, e.g.,~\cite{HeWheWi15}.
\end{remark}

This proposed weak problem formulation is discretized in the following section.

\section{Discretization and a posteriori error estimator}\label{sec_discretization}
In this section, we first discuss the discretization in Section~\ref{discr_form}. Next, 
in Section~\ref{Subsec:Estimator} the residual-type error estimator for the variational inequality in the phase-field model is derived.

\subsection{Discrete formulation}\label{discr_form}
In the discrete setting, at each time step $n = 1,\ldots, N$, we decompose
the polygonal domain $\Omega$ by a (family of) meshes $\mathfrak{M}^n$ consisting
of shape regular parallelograms or triangles $\mathfrak{e}$, such that all meshes share a
common coarse mesh. To allow for local
refinement, in particular of parallelogram elements, we allow for one hanging node per edge at which degrees of freedom
will be eliminated to assert conformity of the discrete
spaces. Further, we assume that the boundary of the domain is resolved by the chosen meshes.

To each mesh, we associate the mesh size function $h^n$, i.e.,
$h^n_{\mathfrak{e}} = \operatorname{diam}{\mathfrak{e}}$ for any
element $\mathfrak{e} \in \mathfrak{M}^n$.
The set of nodes $p$ is given by $\mathfrak{N}$ and we distinguish between the set
$\mathfrak{N}^{\Gamma}$ of nodes at the boundary
and the set of interior nodes $\mathfrak{N}^I$.

Further, for a point $p \in \mathfrak{N}$, we define a patch $\omega_p$ as the interior of the union of all elements sharing the node $p$. We call the union of all sides in the interior of $\omega_p$, not including the boundary of $\omega_p$ skeleton and denote it by $\gamma_p^I$. For boundary nodes, we denote the intersections between $\Gamma$ and $\partial\omega_p$ by  $\gamma_p^{\Gamma}:=\Gamma \cap\partial\omega_p$.
Further, we will make use of $\omega_{\mathfrak{s}}$ which is the union of all elements sharing a side $\mathfrak{s}$.
Later on, we need the definition of the jump term $[\nabla\psi_h]:= \nabla|_{\mathfrak{e}}\psi_h\cdot\boldsymbol{n}_{\mathfrak{e}}- \nabla|_{\tilde{\mathfrak{e}}}\psi_h\cdot\boldsymbol{n}_{\mathfrak{e}}$ where $\mathfrak{e}, \tilde{\mathfrak{e}}$ are neighboring elements and $\boldsymbol{n}_{\mathfrak{e}}$ is the unit outward normal on the common side of the two elements. 

For the discretization, we consider linear finite elements on triangles and bilinear finite elements on parallelograms. 
We abbreviate
\begin{displaymath}
\mathbb{S}_1(\mathfrak{e}):=\left\{\begin{array}{cccc} \mathbb{P}_1(\mathfrak{e}), &\mbox{if}\; &\mathfrak{e} \;\mbox{is a}&\ \mbox{triangle,}\\  \mathbb{Q}_1(\mathfrak{e}), &\mbox{if}\;& \mathfrak{e} \;\mbox{is a}&\ \mbox{parallelogram}. \end{array}\right. 
\end{displaymath}
Thus, discrete solution spaces are given by
\begin{equation*}
  \begin{aligned}
  \boldsymbol{\mathcal{H}}^n_h&:= \{\boldsymbol{v}_h \in \mathcal{C}^0(\overline{\Omega};\mathbb{R}^2)\mid\forall\mathfrak{e}\in\mathfrak{M}^n, \;\boldsymbol{v}_h|_{\mathfrak{e}}\in \mathbb{S}_1(\mathfrak{e})\} \subset \boldsymbol{\mathcal{H}},\\[3pt]
  \boldsymbol{\mathcal{H}}^n_{h,0}&:= \{\boldsymbol{v}_h \in \boldsymbol{\mathcal{H}}^n_h\mid \boldsymbol{v}_h=\boldsymbol{0}\;\mbox{on }\Gamma\}\subset\boldsymbol{\mathcal{H}}_0,\\[3pt]
  \mathcal{H}^n_h&:= \{v_h \in \mathcal{C}^0(\overline{\Omega})\mid\forall\mathfrak{e}\in\mathfrak{M}^n, \;v_h|_{\mathfrak{e}}\in \mathbb{S}_1(\mathfrak{e})\} \subset\mathcal{H}.
  \end{aligned}
\end{equation*}
We define the respective nodal interpolation operators as $I_h^n$, and define the discrete feasible set for the phase-field by 
\begin{equation*}
\mathcal{K}_h^n:=\{\psi_h\in \mathcal{H}^n_h\mid \psi_h(p)\leq (I_h^n\varphi_h^{n-1})(p), \quad \forall p \in \mathfrak{N}\}.
\end{equation*}
The nodal basis functions of the finite element spaces $\mathcal{H}_h^n$ are denoted by $\phi_p$. 

Analogous to Problem~\ref{WeakFormulation}, we define the spatially discretized
time step problem:
\begin{problem}[Discrete formulation in each time step]\label{DiscreteFormulation}
  For $n = 1,\ldots, N$ and $\varphi^{n-1}$ given, find $(\boldsymbol{u}_h^n, \varphi_h^n)\in(\boldsymbol{\mathcal{H}}^n_{h,0}+\{I_h^n\boldsymbol{u}^n_D\})\times\mathcal{K}_h^n$ such that
  \begin{equation}\label{VI_discrete_orig}
\begin{aligned}
\left<g(I_h^n\varphi_h^{n-1})\boldsymbol{\sigma}^+(\boldsymbol{u}_h^n), \boldsymbol{E}_{\mathrm{lin}}(\boldsymbol{w}_h)\right> + \left<\boldsymbol{\sigma}^-(\boldsymbol{u}_h^n) ,\boldsymbol{E}_{\mathrm{lin}}(\boldsymbol{w}_h) \right>& = 0\quad \forall \boldsymbol{w}_h\in \boldsymbol{\mathcal{H}}^n_{h,0},\\[3pt]
\left<(1-\kappa)\varphi_h^n \boldsymbol{\sigma}^+(\boldsymbol{u}_h^n):\boldsymbol{E}_{\mathrm{lin}}(\boldsymbol{u}_h^n),\psi_h -\varphi_h^n\right>  - \frac{G_c}{\epsilon}\left<1-\varphi_h^n, \psi_h -\varphi_h^n \right>&\\[3pt]
 + \epsilon G_c\left<\nabla\varphi_h^n,\nabla(\psi_h -\varphi_h^n)\right>  &\ge 0\quad \forall \psi_h\in\mathcal{K}_h^n.
\end{aligned}
  \end{equation}
\end{problem}
 
\subsection{Residual-type a posteriori estimator}\label{Subsec:Estimator}
In this section, we propose a residual-type a posteriori estimator for
the variational inequality~\eqref{VI_discrete_orig} in Problem~\ref{DiscreteFormulation}. As the structure remains the same for all
time steps, we consider exemplary one time step $n$. We drop the
now superfluous superscript $n$, e.g., $\mathcal{H} :=  \mathcal{H}^n$,  $\mathcal{H}_h := \mathcal{H}^n_h$ and $h_{\mathfrak{e}}:= h^{n}_{\mathfrak{e}}$. 
With the bilinear form 
\begin{equation}\label{EpsBilinearFormDiscrete}
a_{h,\epsilon}(\zeta,\psi) := \left<\left(\frac{G_c}{\epsilon} + \left(1-\kappa\right)\left(\boldsymbol{\sigma}^+\left(\boldsymbol{u}^n_h\right):\boldsymbol{E}_{\mathrm{lin}}\left(\boldsymbol{u}^n_h\right)\right)\right)\zeta,\psi\right>  + G_c\epsilon \left<\nabla \zeta,\nabla \psi \right>,
\end{equation}
we abbreviate the discrete variational inequality in~\eqref{VI_discrete_orig} by
\begin{equation}\label{VI_DiscreteBilinearForm_Discrete}
a_{h,\epsilon}(\varphi_h , \psi_h -\varphi_h ) \ge\left< \frac{G_c}{\epsilon}, \psi_h - \varphi_h\right>\quad \forall \psi_h \in \mathcal{K}^n_h.
\end{equation}
We define the corresponding discrete constraining force density
$\Lambda_h\in \mathcal{H}_h^*$ as
\begin{equation}\label{DiscreteConstrainingForce}
\left< \Lambda_h,\psi_h \right>_{-1,1} := \left<\frac{G_c}{\epsilon},\psi_h\right>-a_{h,\epsilon}(\varphi_h, \psi_h)\quad\forall \psi_h\in\mathcal{H}_{h}.
\end{equation}
We note that the discrete constraining force $\Lambda_h$ equals the linear residual $\mathcal{R}_h^{lin}$ for the corresponding unconstrained diffusion-reaction equations.

Further, with the bilinear form~\eqref{EpsBilinearFormDiscrete} and
the admissible set $\mathcal{K}(I^n_h(\varphi^{n-1}_h)):=\{\psi\in
\mathcal{H}\mid \psi\leq I^n_h(\varphi_h^{n-1})\}$, we approximate the
variational inequality in~\eqref{VI_orig} by introducing an auxiliary variable
$\hat{\varphi} \in \mathcal{K}(I^n_h(\varphi^{n-1}_h))$ solving
\begin{equation}\label{VI_DiscreteBilinearForm}
a_{h,\epsilon}(\hat{\varphi} , \psi -\hat{\varphi} ) \ge \left< \frac{G_c}{\epsilon},\psi - \hat{\varphi}  \right>\quad \forall \psi \in \mathcal{K}(I_h(\varphi^{n-1}_h)).
\end{equation}
\begin{remark}
As the bilinear form $a_{h,\epsilon}(\cdot,\cdot)$ depends on the approximation $\boldsymbol{u}^n_h$ of $\boldsymbol{u}^n$ and the constraints depend on the approximation $I^n_h(\varphi_h^{n-1})$  of $\varphi^{n-1}$, the solution $\hat{\varphi}$ of~\eqref{VI_DiscreteBilinearForm} 
is an approximation to the solution $\varphi$ of~\eqref{VI_orig}. 
\end{remark}

It exists a distribution $\hat{\Lambda}\in\mathcal{H}^{\ast}$, called constraining force density, which turns the variational inequality~\eqref{VI_DiscreteBilinearForm} in an equation
\begin{equation*}
\left<\hat{\Lambda} ,\psi\right>_{-1,1}:=\left<\frac{G_c}{\epsilon}, \psi \right> -a_{h,\epsilon}(\hat{\varphi}, \psi )\quad \forall \psi\in \mathcal{H}.
\end{equation*}

As the work ~\cite{Veeser:2001} reveals that sharp a posteriori estimators for variational inequalities can be derived by involving the error in the constraining forces, we measure the error in the solution of~\eqref{VI_DiscreteBilinearForm} as well as in the constraining forces. Similar to~\cite{Walloth:2018b}, we measure the error of the solution of~\eqref{VI_DiscreteBilinearForm} in the energy norm
\begin{equation*}
\|\cdot\|_{\epsilon}:= \left\{G_c\epsilon\|\nabla (\cdot)\|^2 + \left\|\left(\frac{G_c}{\epsilon}+\left(1-\kappa\right)\boldsymbol{\sigma}^+\left(\boldsymbol{u}^n_h\right):\boldsymbol{E}_{\mathrm{lin}}\left(\boldsymbol{u}^n_h\right)\right)^{\frac{1}{2}}\left(\cdot\right)\right\|^2 \right\}^{\frac{1}{2}},
\end{equation*}
corresponding to the bilinear form $a_{h,\epsilon}(\cdot,\cdot)$. The
error in the constraining forces is measured in the corresponding dual
norm
\[
  \|\cdot\|_{\ast,\epsilon}:= \frac{\mathrm{sup}_{\psi\in
      \mathcal{H}^1}\left<\cdot, \psi \right>}{\|\psi\|_{\epsilon}}.
\]
In order to compare the continuous and discrete constraining forces,
we cannot use $\Lambda_h$, as by
definition~\eqref{DiscreteConstrainingForce} it is a functional on the
space of discrete functions and there is no unique extension to $\mathcal{H}^*$. Therefore, we need to choose an extension which is a discrete counterpart of $\hat{\Lambda}$ as functional on $H^1$, reflecting the properties of the constraining force and depending on the discrete solution and given data. We call this extension quasi-discrete constraining force and denote it by $\tilde{\Lambda}_h$.

In~\cite{Veeser:2001}, an extension to a functional on $H^1$ was proposed by means of lumping 
$\sum_{p\in\mathfrak{N}^C}s_p\phi_p$,
where 
\begin{align}\label{NodeValueLumpedForce}
s_p = \frac{\left<\Lambda_h,\phi_p\right>_{-1,1}}{\int_{\omega_p}\phi_p\ \mathrm{d}x}\ge 0
\end{align}
are the node values of the lumped discrete constraining force. 
This approach has been extended and applied to different obstacle and contact problems in, e.g.,~\cite{NochettoSiebertVeeser:2005, MoonNochettoPetersdorffZhang:2007, KrauseVeeserWalloth:2015, GudiPorwal:2014, GudiPorwal:2016,Walloth:2019,Walloth:2018a,Walloth:2018b}.

Following the works~\cite{MoonNochettoPetersdorffZhang:2007, Walloth:2018b}, we distinguish between full-contact nodes
$p\in\mathfrak{N}^{fC}$ and semi-contact nodes
$p\in\mathfrak{N}^{sC}$. Full-contact nodes are those nodes for which
the solution is fixed to the obstacle $\varphi_h=I_h^n\varphi^{n-1}_h$
on $\omega_p$ and the sign condition
\[
  \left<\mathcal{R}^{lin}_h,\varphi \right>_{-1,1,\omega_p}\ge 0 \quad
  \forall \varphi\ge 0\in H^1_0(\omega_p)
\]
is fulfilled. The latter condition means that the solution is locally not improvable, see the explanation in~\cite{MoonNochettoPetersdorffZhang:2007}. Semi-contact nodes are those nodes for which $\varphi_h(p)=I_h^{n}\varphi^{n-1}_h(p)$ holds but not the conditions of full-contact. 
Based on this classification, we define the quasi-discrete constraining force as
\begin{align}\label{QuasiDiscreteConstrainingForce}
&\left< \tilde{\Lambda}_h,\psi \right>_{-1,1} := \sum_{p\in\mathfrak{N}^{sC}}\left< \tilde{\Lambda}_h^p,\psi\phi_p \right>_{-1,1} + \sum_{p\in\mathfrak{N}^{fC}}\left< \tilde{\Lambda}_h^p,\psi\phi_p \right>_{-1,1}.
\end{align}
For the definition of the local contributions, we shorten the element residual with
\begin{equation*}
r(\varphi_h) :=\frac{G_c}{\epsilon} + G_c\epsilon\Delta \varphi_h -  \frac{G_c}{\epsilon}\varphi_h- (1-\kappa)(\boldsymbol{\sigma}^+(\boldsymbol{u}^n_h):\boldsymbol{E}_{\mathrm{lin}}(\boldsymbol{u}^n_h))\varphi_h.
\end{equation*}
For semi-contact nodes, we consider the following local contribution in~\eqref{QuasiDiscreteConstrainingForce}
\begin{align*}
\left< \tilde{\Lambda}_h^p,\psi\phi_p \right>_{-1,1}=\int_{\gamma_p^I}G_c\epsilon[\nabla \varphi_h]c_p(\psi)\phi_p\ \mathrm{d}x  
+ \int_{\omega_p} r(\varphi_h)c_p(\psi)\phi_p \ \mathrm{d}x,
\end{align*} 
with $c_p(\psi)= \frac{\int_{\tilde{\omega}_p}\psi\phi_p\ \mathrm{d}x}{\int_{\tilde{\omega}_p}\phi_p\ \mathrm{d}x}$, where $\tilde{\omega}_p$ is the patch around $p$ with respect to two uniform red-refinements. 
For full-contact nodes we define the following local contribution in~\eqref{QuasiDiscreteConstrainingForce} as
\begin{align*}
\left< \tilde{\Lambda}_h^p,\psi\phi_p \right>_{-1,1} :=\int_{\gamma_p^I}G_c\epsilon^2[\nabla \varphi_h]\psi\phi_p\ \mathrm{d}x   + \int_{\omega_p} r(\varphi_h)\psi\phi_p\ \mathrm{d}x. 
\end{align*}
With these definitions, we define the error measure 
\begin{equation}\label{ErrorMeasure}
\|\hat{\varphi}-\varphi_h\|_{\epsilon} + \|\hat{\Lambda}-\tilde{\Lambda}_h\|_{\ast,\epsilon}.
\end{equation}
In order to state the error estimator 
\begin{equation}\label{Def_Estimator}
\eta := \sum_{k=1}^4\eta_{k}
\end{equation}
for the error measure~\eqref{ErrorMeasure}, we define for each node $p$
\begin{equation*}
\alpha_p:= \mathrm{min}_{x\in\omega_p}\left\{\frac{G_c}{\epsilon}+(1-\kappa) (\boldsymbol{\sigma}^+(\boldsymbol{u}^n_h):\boldsymbol{E}_{\mathrm{lin}}(\boldsymbol{u}^n_h))\right\},
\end{equation*}
and $h_p:=\mathrm{diam}(\omega_p)$.
The local contributions are 
\begin{align*}
\eta_{1}:=&\left(\sum_{p\in\mathfrak{N}\backslash\mathfrak{N}^{fC}}\eta^2_{1,p}\right)^{\frac{1}{2}}, & \eta_{1,p}:=&\mathrm{min}\left\{\frac{h_p}{\sqrt{G_c\epsilon}},\alpha_p^{-\frac{1}{2}}\right\}\|r(\varphi_h)\|_{\omega_p},\\[3pt] 
\eta_{2}:=&\left(\sum_{p\in\mathfrak{N}^I\backslash\mathfrak{N}^{fC}}\eta^2_{2,p}\right)^{\frac{1}{2}}, &\eta_{2,p}:=& \mathrm{min}\left\{\frac{h_p}{\sqrt{G_c\epsilon}},\alpha_p^{-\frac{1}{2}}\right\}^{\frac{1}{2}}(G_c\epsilon)^{-\frac{1}{4}}\|G_c\epsilon[\nabla \varphi_h]\|_{\gamma_p^I},\\[3pt] 
\eta_{3}:=&\left(\sum_{p\in\mathfrak{N}^{\Gamma}\backslash\mathfrak{N}^{fC}}\eta^2_{3,p}\right)^{\frac{1}{2}}, &\eta_{3,p}:=& \mathrm{min}\left\{\frac{h_p}{\sqrt{G_c\epsilon}},\alpha_p^{-\frac{1}{2}}\right\}^{\frac{1}{2}}(G_c\epsilon)^{-\frac{1}{4}}\|G_c\epsilon \nabla \varphi_h \|_{\gamma_p^{\Gamma}}, \\[3pt] 
\eta_{4}:=&\left(\sum_{p\in\mathfrak{N}^{sC}}\eta^2_{4,p}\right)^{\frac{1}{2}}, &\eta_{4,p}:=&\left(s_p \int_{\tilde{\omega}_p}(I_h^n\varphi^{n-1}_h-\varphi_h)\phi_p\ \mathrm{d}x\right)^{\frac{1}{2}}.
\end{align*}
We emphasize that the estimator contributions related to the
constraints are localized to the area of semi-contact
  and no contributions arise from full-contact nodes. 

Under the assumption  that $\varphi_h$ is a linear finite
element function, we will prove in a forthcoming publication that $\eta$ constitutes global upper and local lower bounds of~\eqref{ErrorMeasure}.
The estimator $\eta$ provides  a robust upper bound where robust means that
the constant in the bound does not depend on
$\epsilon$ such that the validity of the estimator holds for arbitrary choices of $\epsilon$.
\begin{theorem}{\bf Reliability of the error estimator}
\label{Theorem:UpperBound}\\
The error estimator $\eta$~\eqref{Def_Estimator} provides a robust upper bound of the error measure~\eqref{ErrorMeasure}:
\begin{equation*}
\|\hat{\varphi}-\varphi_h\|_{\epsilon} + \|\hat{\Lambda}-\tilde{\Lambda}_h\|_{\ast,\epsilon}\lesssim \eta.
\end{equation*}
\end{theorem} 

The local lower bounds are summarized in the following Theorems.
\begin{theorem}{\bf Local lower bounds by $\eta_{1,p},\eta_{2,p},\eta_{3,p}$}\label{Theorem:LowerBound}\\
The error estimator contributions $\eta_{k,p}$, $k=1,2,3$ constitute the following robust local lower bounds
\begin{equation*}
\eta_{k,p}\lesssim\|\hat{\varphi}-\varphi_h\|_{\epsilon,\omega_p} + \|\hat{\Lambda}-\tilde{\Lambda}_h\|_{\ast,\epsilon,\omega_p}.
\end{equation*}
\end{theorem}
\begin{theorem}{\bf Local lower bound by $\eta_{4,p}$}\label{Theorem:LowerBound2}\\
For nodes $p\in\mathfrak{N}^{sC}$ with $\frac{h_p}{\sqrt{G_c\epsilon}} \leq \alpha_p^{-\frac{1}{2}}$ we have the robust local lower bound
\begin{equation}\label{LowerBound4_Robust}
\begin{split}
\eta_{4,p} \lesssim& \|\hat{\varphi}-\varphi_h\|_{\epsilon,\omega_p} + \|\hat{\Lambda}-\tilde{\Lambda}_h\|_{\ast,\epsilon,\omega_p}  \\[3pt]
&  + \mathrm{min}\left\{\frac{h_p}{\sqrt{G_c\epsilon}},\alpha_p^{-\frac{1}{2}}\right\}^{\frac{1}{2}}(G_c\epsilon)^{-\frac{1}{4}}\|G_c\epsilon[\nabla (I_h^n\varphi^{n-1}_h-\varphi_h)]^I\|_{\gamma_p}.
\end{split}
\end{equation}
Otherwise, for nodes $p\in\mathfrak{N}^{sC}$ with $\alpha_p^{-\frac{1}{2}}< \frac{h_p}{\sqrt{G_c\epsilon}}$
we have the local lower bound
\begin{equation}\label{LowerBound4_NonRobust}
\begin{split}
\eta_{4,p} \lesssim &\|\hat{\varphi}-\varphi_h\|_{\epsilon,\omega_p} + \|\hat{\Lambda}-\tilde{\Lambda}_h\|_{\ast,\epsilon,\omega_p} \\[3pt]
& + \alpha_p^{\frac{3}{4}}(G_c\epsilon)^{-\frac{3}{4}} \mathrm{min}\left\{\frac{h_p}{\sqrt{G_c\epsilon}},\alpha_p^{-\frac{1}{2}}\right\}^{\frac{1}{2}}(G_c\epsilon)^{-\frac{1}{4}}\|G_c\epsilon[\nabla (I_h^n\varphi^{n-1}_h-\varphi_h)]^I\|_{\gamma_p}.
\end{split}
\end{equation}
\end{theorem}
\begin{remark}
We note that the additional term in the bound~\eqref{LowerBound4_Robust} only occurs for $p\in\mathfrak{N}^{sC}$ and is of the same order as the other estimator contributions. 
In the application, we expect the semi-contact zone to be well resolved, especially with respect to $\epsilon$ such that $\frac{h_p}{\sqrt{G_c\epsilon}} \leq \alpha_p^{-\frac{1}{2}}$ 
after a finite number of adaptive refinement steps such that the local lower bound is robust everywhere. 
 \end{remark}

\begin{remark}
  In this work, we focus on the novel estimator for the phase-field
  inequality to obtain a  good resolution of the fracture growth.
  We therefore only provide a residual-type a posteriori estimator for
  the variational inequality in Problem~\ref{DiscreteFormulation}.
  For the equation in the coupled system of
  Problem~\ref{DiscreteFormulation},
  a standard a posteriori estimator~\cite{Verfuerth_1998a} could be applied. 
\end{remark}

\section{Solver and refinement strategy}\label{sec:solver}
The numerical solution proceeds from Problem~\ref{DiscreteFormulation}.
Concerning the robustness, efficiency and the accuracy of the coupling terms, we made good experiences
treating the phase-field system in a monolithic fashion, e.g.,~\cite{Wi17_SISC,Wi17_CMAME}. At first, the handling of the 
crack irreversibility constraint is clarified. In the following section, the spatial
discretization and the overall solution method are explained. 


\subsection{Solution algorithms}
In order to compute a discrete approximation of the solution of the quasi-static fracture phase-field model, we use a semi-smooth Newton method~\cite{Hintermueller_Ito_Kunisch_2003} implemented in~\cite{DOpElib}. The semi-smooth Newton method is based on the following complementarity system which is equivalent to the variational inequality system of Problem~\ref{WeakFormulation}. 

It is easy to see, that by introducing a Lagrange multiplier
$\Lambda^n \in
\mathcal{H}^*$~\cite{ito2008lagrange,rockafellar1993lagrange,vexler2008adaptive}
the weak formulation~\eqref{VI_orig}, in the continuous setting, is
equivalent to the complementarity system: given $n=1,\ldots, N$
and $\varphi^{n-1}$
to find $(\boldsymbol{u}^n, \varphi^n,
\Lambda^n)\in(\boldsymbol{\mathcal{H}}_0+\{\boldsymbol{u}^n_D\})\times\mathcal{K}(\varphi^{n-1})\times
\mathcal{H}^*$ satisfying
\begin{equation*}
  \begin{aligned}
    \left<g(\varphi^{n-1})\boldsymbol{\sigma}^+(\boldsymbol{u}^n), \boldsymbol{E}_{\mathrm{lin}}(\boldsymbol{w})\right> + \left<\boldsymbol{\sigma}^-(\boldsymbol{u}^n) ,\boldsymbol{E}_{\mathrm{lin}}(\boldsymbol{w}) \right>& = 0\ &&\forall \boldsymbol{w}\in \boldsymbol{\mathcal{H}}_0,\\[3pt]
    \left<(1-\kappa)\varphi^n \boldsymbol{\sigma}^+(\boldsymbol{u}^n):\boldsymbol{E}_{\mathrm{lin}}(\boldsymbol{u}^n),\psi\right>  - \frac{G_c}{\epsilon}\left<1-\varphi^n, \psi \right>&\\[3pt]
    + \epsilon G_c\left<\nabla\varphi^n,\nabla\psi\right> + \langle \Lambda^n,\psi\rangle_{-1,1} &= 0\ &&\forall \psi\in\mathcal{H},\\[3pt]
    \langle \Lambda^n,\psi\rangle_{-1,1} &\ge 0,\ &&\forall \psi \in
    \mathcal{H}, \psi \ge 0,\\[3pt]
    \varphi^n &\le \varphi^{n-1},\\
    \langle \Lambda^n,\varphi^{n-1}-\varphi^n\rangle_{-1,1} &= 0.
  \end{aligned}
\end{equation*}

To obtain the complementarity formulation for the
discretization~\eqref{VI_discrete_orig}, we define
$(\mathcal{H}^n_h)^*$ by a dual basis of $\mathcal{H}^n_h$; i.e.,
we let $(\mathcal{H}^n_h)^* = \operatorname{span}\{\phi^*_p\mid p \in \mathfrak{N}\}$, where
\[
\langle \phi_p^*, \phi_q\rangle = \delta_{pq}
\]
for the nodal basis $\phi_q$ of $\mathcal{H}^n_h$.
Then, we define $\Lambda_h^n = \sum_{p \in \mathfrak{N}}
(\Lambda_h^n)_p \phi_p^*$ by setting
\[
\langle \Lambda_h^n,\phi_p\rangle =
-\left<(1-\kappa)\varphi_h^n \boldsymbol{\sigma}^+(\boldsymbol{u}_h^n):\boldsymbol{E}_{\mathrm{lin}}(\boldsymbol{u}_h^n),\phi_p\right>
 + \frac{G_c}{\epsilon}\left<1-\varphi_h^n, \phi_p \right>
 - \epsilon G_c\left<\nabla\varphi_h^n,\nabla\phi_p\right> 
 \]
 for all $p \in \mathfrak{N}$. This immediately gives the analogous discrete
 complementarity system
\begin{equation}\label{DiscretComplementarity}
  \begin{aligned}
    \left<g(I_h^n\varphi_h^{n-1})\boldsymbol{\sigma}^+(\boldsymbol{u}_h^n), \boldsymbol{E}_{\mathrm{lin}}(\boldsymbol{w}_h)\right> + \left<\boldsymbol{\sigma}^-(\boldsymbol{u}_h^n) ,\boldsymbol{E}_{\mathrm{lin}}(\boldsymbol{w}_h) \right>& = 0\ &&\forall \boldsymbol{w}_h\in \boldsymbol{\mathcal{H}}^n_{h,0}\\[3pt]
    \left<(1-\kappa)\varphi_h^n \boldsymbol{\sigma}^+(\boldsymbol{u}_h^n):\boldsymbol{E}_{\mathrm{lin}}(\boldsymbol{u}_h^n),\psi_h\right>  - \frac{G_c}{\epsilon}\left<1-\varphi_h^n, \psi_h \right>&\\[3pt]
    + \epsilon G_c\left<\nabla\varphi_h^n,\nabla\psi_h\right> + \langle \Lambda_h^n,\psi_h\rangle &= 0\quad &&\forall \psi_h\in\mathcal{H}^n_h,\\[3pt]
    \langle \Lambda_h^n,\phi_p\rangle &\ge 0,\ &&\forall p \in \mathfrak{N},\\[3pt]
    \varphi^n(p) &\le \varphi^{n-1}(p),\ &&\forall p \in \mathfrak{N}, \\[3pt]
    \langle \Lambda_h^n,I_h^n\varphi^{n-1}-\varphi^n\rangle &= 0.
  \end{aligned}
\end{equation}

\begin{remark}
We recall that the discrete variational
inequality~\eqref{VI_DiscreteBilinearForm_Discrete} for which we
derived the a posteriori estimator in Section~\ref{Subsec:Estimator}
is part of the discrete phase-field model
(Problem~\ref{DiscreteFormulation}). Due to the equivalence of
Problem~\ref{DiscreteFormulation} and the complementarity
system~\eqref{DiscretComplementarity} it holds that $\left<
  \Lambda_h,\phi_p \right>_{-1,1}$ in~\eqref{DiscretComplementarity}
and~\eqref{DiscreteConstrainingForce}  are the same. The only
difference is that now we have chosen a discrete basis for
$\Lambda_h$ such that the node value $(\Lambda_h^n)_p$ equals
$s_p\cdot \int_{\omega_p}\phi_p\ \mathrm{d}x$ in ~\eqref{NodeValueLumpedForce}.
\end{remark}

To apply a semi-smooth Newton method, we notice, that by choice of the
basis the complementarity conditions are equivalent to the following
complementarity system for the coefficients $\varphi^n(p)$ and $(\Lambda_h^n)_p$:
\begin{align*}
\begin{aligned}
     (\Lambda_h^n)_p = \langle \Lambda_h^n,\phi_p\rangle &\ge 0,\ &&\forall p \in \mathfrak{N},\\[3pt]
    \varphi^n(p) &\le \varphi^{n-1}(p),\ &&\forall p \in \mathfrak{N}, \\
    \sum_{p \in \mathfrak{N}} (\Lambda_h^n)_p(I_h^n\varphi^{n-1}-\varphi^n)(p) = \langle \Lambda_h^n,I_h^n\varphi^{n-1}-\varphi^n\rangle &= 0.
  \end{aligned}
  \end{align*}
  Using the complementarity function $(x,y) \mapsto x - \max(0,x+cy)$,
  for arbitrary $c > 0$, we can equivalently express~\eqref{DiscretComplementarity} as the system
  \begin{equation}\label{NewtonSystem}
    \begin{aligned}
      \left<g(I_h\varphi_h^{n-1})\boldsymbol{\sigma}^+(\boldsymbol{u}_h^n), \boldsymbol{E}_{\mathrm{lin}}(\boldsymbol{w}_h)\right> + \left<\boldsymbol{\sigma}^-(\boldsymbol{u}_h^n) ,\boldsymbol{E}_{\mathrm{lin}}(\boldsymbol{w}_h) \right>& = 0\ &&\forall \boldsymbol{w}_h\in \boldsymbol{\mathcal{H}}^n_{h,0}\\[3pt]
      \left<(1-\kappa)\varphi_h^n \boldsymbol{\sigma}^+(\boldsymbol{u}_h^n):\boldsymbol{E}_{\mathrm{lin}}(\boldsymbol{u}_h^n),\psi_h\right>  - \frac{G_c}{\epsilon}\left<1-\varphi_h^n, \psi_h \right>&\\[3pt]
      + \epsilon G_c\left<\nabla\varphi_h^n,\nabla\psi_h\right> + \langle \Lambda_h^n,\psi_h\rangle_{-1,1} &= 0\ &&\forall \psi_h\in\mathcal{H}^n_h,\\[3pt]
      (\Lambda_h^n)_p - \max(0,(\Lambda_h^n)_p + c((I_h^n\varphi^{n-1}-\varphi^n)(p))) &= 0,\ &&\forall p \in \mathfrak{N}.
    \end{aligned}
  \end{equation}
  Now, a time step $n$ can be reformulated,
  using~\eqref{NewtonSystem}, as given $\varphi_h^{n-1}$:\\
  Find $U_h^n:= U_h = (u_h,\varphi_h,\Lambda_h) \in \mathcal{U}_h + \{\boldsymbol{u}_D^n\}:=
  ((\boldsymbol{\mathcal{H}}_{h,0} + \{\boldsymbol{u}_D^n\}) \times
  \mathcal{H}_h\times \mathcal{H}_h^*)$, such that
  \[
   A(U_h) = 0
 \]
 abbreviating~\eqref{NewtonSystem}.
  To solve this non-linear equation, we formulate 
  a residual-based Newton scheme utilizing that the $\max$ operator is semi-smooth in finite
  dimensions. 
  
  To measure the residuals and monitoring functions, we use 
  the discrete norm $\|\cdot\| :=\|\cdot\|_{l^\infty}$ measuring the
  maximal absolute value of the coordinate vectors.
  At a given time instance $t_n$, we shall find the loading step solution $U_h^n$ using:
  \begin{Algorithm}[Residual-based Newton's method]
    \label{algo_residual_based_Newton}
    Choose $\rho \in (0,1)$ and an initial Newton guess $U^{(0)} \in \mathcal{U}_h +
    \{\boldsymbol{u}_D^n\}$. Typically, $U^{(0)} =  I_h^nU_h^{n-1}$ up
    to a correction for the Dirichlet-values.
    
    For the iteration steps $k=0,1,2,3,\ldots$:
    \begin{enumerate}
    \item Stop, if converged and return $U^n_h = U^{(k)}$.
    \item Find $\delta U^{(k)} \in \mathcal{U}_h$ such that
      \begin{align*}
        A'(U^k)(\delta U^k) &= - A(U^k).
      \end{align*}
    \item Find a maximal $s_k \in \{\rho^l\mid l \in 0,\ldots,1 \}$
      such that 
      \begin{equation*}
        \| A(U^{k}+s_k \delta U^{(k)})\| < \|A(U^k)\|.
      \end{equation*}
    \item Set 
      \[
        U^{(k+1)} = U^{(k)} + s_k \delta U^{(k)}.
      \]
    \end{enumerate}
\end{Algorithm}
For computational reasons, several of the steps need to
  be slightly modified, e.g., by only reassembling the matrix if the
  convergence rate is too slow, and for simplicity the linear systems
  in the Newton algorithm are solved by UMFPACK~\cite{davis2004algorithm}.
  For details, we refer to our
  implementation in~\cite{dope,DOpElib}.


\subsection{Refinement strategy}\label{sec_refine}
The used adaptive solution strategy is given in this section.
The mesh adaptation is realized using extracted local error indicators from the 
a posteriori error estimator~\eqref{Def_Estimator} on the given meshes
over all time steps.

This information is used to adapt the mesh using the following 
strategy:
\begin{Algorithm}\label{alg:refinement}
  Given a time discretization $t_0 < \ldots< t_N$, and an initial
  mesh $\mathfrak{M}^n = \mathfrak{M}$ for all $n = 0, \ldots, N$.
 \begin{enumerate}
 \item Set $\varphi^0_h = I_h^0 \varphi^0$ and
   solve the discrete complementarity system~\eqref{NewtonSystem} to 
   obtain the discrete solutions $\boldsymbol{u}^n_h,\varphi^n_h,\Lambda_h^n$ for all $n=1,\ldots,N$.
 \item Evaluate the error estimator~\eqref{Def_Estimator} in order to
   obtain $\eta^n$ for each time step.
 \item Stop, if $\sum_{n=1}^N(\eta^n)^2$ and $\|I_h^n
   \varphi^{n-1}-\varphi^{n-1}\|$ are small enough for all $n = 1,\ldots,
   N$.
\item For each $n=1,\ldots N$, mark elements in $\mathfrak{M}^n$ based on $\eta^n$ according to an
  optimization strategy, as implemented in deal.II~\cite{dealII85}.
  This strategy allows to flag certain cells to reach a grid that is optimal with respect to an objective function that tries to balance reducing the error and increasing the numerical cost. 
  More details on this approach can be found in~\cite{richter2005parallel}.
\item Refine the meshes according to the marking, and satisfaction of
  the constraints on hanging nodes.
\item Repeat from step 1.
\end{enumerate}
\end{Algorithm}


\section{Numerical tests}\label{sec_tests}

In this section, we study the quality of the error estimator proposed in
Section~\ref{Subsec:Estimator}. 
For that, we use two crack propagation settings in pure elasticity
regimes (each with three studies resulting in a total of six scenarios).
The following questions and aspects are addressed:
\begin{itemize}
 \item Does the error estimator allow to resolve the mushy zone around the crack sufficiently? (Both Studies in Section~\ref{study_one} and~\ref{study_two})
 \item How does adaptive mesh refinement performs in comparison to uniform
   mesh refinement in terms of the convergence of the, so-called,
   load-displacement curves? (Study $1$ in Section~\ref{study_one})
 \item Investigating the $\epsilon-h$ relationship (Study $2$ in Section~\ref{study_two})
 \item Observing the error indicators and the corresponding adaptively
  refined meshes (Study $3$ in Section~\ref{study_three}).
\end{itemize}
The implementation is done in the open-source package Differential
  Equations and Optimization Environment library
  (DOpElib)~\cite{dope,DOpElib} using the finite elements from deal.II~\cite{dealII85}.

In the simulations, we declare a strip of size $4 h_{\text{start}}$
below the top boundary to be ignored by the estimator. Thus, we avoid that the error
estimator resolves the singularity due to the non-matching boundary
conditions which allow for a fracture to form directly below the top-boundary.  

\subsection{Configurations}\label{config}
The two numerical configurations are set according to~\cite{miehe2010phase}: the single edge
notched shear test and the single edge notched tension test,
the boundary-values are selected according to~\cite{Wi17_CMAME}.
Both tests were used by several groups with similar settings and it is well known that under constant tension
the crack grows in a straight line, while under constant shear forces the
crack grows in a curve towards a 
corner~\cite{borden2012phase,miehe2010phase,HeWheWi15,AmGeraLoren15}.


\subsubsection{The single edge notched shear test}\label{config_shear}
The geometry and the material parameters of the single edge notched shear test are adopted from~\cite{miehe2010phase} and displayed in Figure~\ref{shear_geo}. 
Here, the domain $\Omega$ is a two-dimensional square of $10\ \si{mm}$ length
with a given crack (called slit) on the right side at $5\ \si{mm}$ tending to the midpoint of the square. 
On the bottom boundary the square is fixed, on the top boundary a given force in $x$-direction pulls to the left. We follow the boundary conditions described in~\cite{Wi17_CMAME}.

\begin{figure}[htbp!]
\begin{minipage}[b]{0.48\textwidth}
\centering
 \begin{tikzpicture}[xscale=0.85,yscale=0.85]
\draw[fill=gray!30] (0,0) -- (0,5) -- (5,5) -- (5,0) -- (0,0);
\draw[line width =2pt, blue] (0,0) -- (0,5) -- (5,5) -- (5,0) -- (0,0);
\draw[red, thick] (2.5,2.5) -- (5,2.5);
\draw[line width =2pt, blue] (2.5,2.45) -- (5,2.45);
\node at (3.725,2.15) {slit};
\draw (5,0) -- (5,5);
\draw[line width =2pt, blue] (5,5) -- (0,5);
\draw[-] (0,5) -- (0,0);
\draw[<->] (-0.25,0) -- (-0.25,5);
\draw[->,blue] (0.7,5.5) -- (0.1,5.5);
\node at (2.5,5.3) {$\Gamma_{\text{top}}$};
\node at (0.4,5.25) {$u_x$};
\node at (-0.9,2.5) {$10\ \si{mm}$};
\draw[<->] (0,-0.5) -- (5,-0.5);
\node at (5.15,-0.5) {$x$};
\node at (-0.5,5.15) {$y$};
\node at (2.5,-0.75) {$10\ \si{mm}$};

\draw[<->] (5.25,0) -- (5.25,2.5);
\node at (5.77,1.25) {$5\ \si{mm}$};

\draw (0.1,0) -- (-0.05,-0.2);
\draw (0.3,0) -- (0.15,-0.2);
\draw (0.5,0) -- (0.35,-0.2);
\draw (0.7,0) -- (0.55,-0.2);
\draw (0.9,0) -- (0.75,-0.2);
\draw (1.1,0) -- (0.95,-0.2);
\draw (1.3,0) -- (1.15,-0.2);
\draw (1.5,0) -- (1.35,-0.2);
\draw (1.7,0) -- (1.55,-0.2);
\draw (1.9,0) -- (1.75,-0.2);
\draw (2.1,0) -- (1.95,-0.2);
\draw (2.3,0) -- (2.15,-0.2);
\draw (2.5,0) -- (2.35,-0.2);
\draw (2.7,0) -- (2.55,-0.2);
\draw (2.9,0) -- (2.75,-0.2);
\draw (3.1,0) -- (2.95,-0.2);
\draw (3.3,0) -- (3.15,-0.2);
\draw (3.5,0) -- (3.35,-0.2);
\draw (3.7,0) -- (3.55,-0.2);
\draw (3.9,0) -- (3.75,-0.2);
\draw (4.1,0) -- (3.95,-0.2);
\draw (4.3,0) -- (4.15,-0.2);
\draw (4.5,0) -- (4.35,-0.2);
\draw (4.7,0) -- (4.55,-0.2);
\draw (4.9,0) -- (4.75,-0.2);
 \end{tikzpicture}
 \caption{\textbf{(Single edge notched shear test)}\\ Geometry and boundary conditions. On the left and right side and the lower part of the slit, the boundary condition in $y$-direction is $u_y = 0\ \si{mm}$ and traction-free in $x$-direction. 
 On the bottom boundary it holds $u_x = u_y = 0\ \si{mm}$. On the top boundary, it holds
 $u_y = 0\ \si{mm}$ and in $x$-direction we determine a time-dependent non-homogeneous Dirichlet condition: $u_x = t \cdot 1\ \si{mm/s}$ with $0 \leq t \leq 0.0125\ \si{s}$ with a time step size $\delta t > 0$.}\label{shear_geo}
\end{minipage}
\hfill
\begin{minipage}[b]{0.48\textwidth}
\centering
 \begin{tikzpicture}[xscale=0.85,yscale=0.85]
\draw[fill=gray!30] (0,0) -- (0,5) -- (5,5) -- (5,0) -- (0,0);
\draw[line width =2pt, blue] (0,0) -- (5,0);
\draw[red, thick] (2.5,2.5) -- (5,2.5);
\node at (3.725,2.2) {slit};
\draw (5,0) -- (5,5);
\draw[line width =2pt, blue] (5,5) -- (0,5);
\draw[-] (0,5) -- (0,0);
\draw[<->] (-0.25,0) -- (-0.25,5);
\draw[->,blue] (0.5,5.5) -- (0.5,5.9);
\node at (2.5,5.3) {$\Gamma_{\text{top}}$};
\node at (0.4,5.27) {$u_y$};
\node at (-0.9,2.5) {$10\ \si{mm}$};
\draw[<->] (0,-0.5) -- (5,-0.5);
\node at (5.15,-0.5) {$x$};
\node at (-0.5,5.15) {$y$};
\node at (2.5,-0.75) {$10\ \si{mm}$};

\draw[<->] (5.25,0) -- (5.25,2.5);
\node at (5.77,1.25) {$5\ \si{mm}$};

\draw (0.1,0) -- (-0.05,-0.2);
\draw (0.3,0) -- (0.15,-0.2);
\draw (0.5,0) -- (0.35,-0.2);
\draw (0.7,0) -- (0.55,-0.2);
\draw (0.9,0) -- (0.75,-0.2);
\draw (1.1,0) -- (0.95,-0.2);
\draw (1.3,0) -- (1.15,-0.2);
\draw (1.5,0) -- (1.35,-0.2);
\draw (1.7,0) -- (1.55,-0.2);
\draw (1.9,0) -- (1.75,-0.2);
\draw (2.1,0) -- (1.95,-0.2);
\draw (2.3,0) -- (2.15,-0.2);
\draw (2.5,0) -- (2.35,-0.2);
\draw (2.7,0) -- (2.55,-0.2);
\draw (2.9,0) -- (2.75,-0.2);
\draw (3.1,0) -- (2.95,-0.2);
\draw (3.3,0) -- (3.15,-0.2);
\draw (3.5,0) -- (3.35,-0.2);
\draw (3.7,0) -- (3.55,-0.2);
\draw (3.9,0) -- (3.75,-0.2);
\draw (4.1,0) -- (3.95,-0.2);
\draw (4.3,0) -- (4.15,-0.2);
\draw (4.5,0) -- (4.35,-0.2);
\draw (4.7,0) -- (4.55,-0.2);
\draw (4.9,0) -- (4.75,-0.2);
 \end{tikzpicture}
 \caption{\textbf{(Single edge notched tension test)}\\ Geometry and boundary conditions. On the left and right side,
 the boundaries are traction-free (homogeneous Neumann condition).  On the bottom boundary it holds $u_y = 0\ \si{mm}$. On the top boundary, it holds
 $u_x = 0\ \si{mm}$ and in $y$-direction we determine a time-dependent non-homogeneous Dirichlet condition: $u_y = t \cdot 1\ \si{mm/s}$ with $0 \leq t \leq 0.00676\ \si{s}$ with a time step size $\delta t > 0$.}\label{tension_geo}
\end{minipage}
\end{figure}

The material and model parameters are given as follows:
the Lam\'{e} coefficients
are given as $\lambda=121.15\ \si{kN/mm^2}$ and $\mu=80.77\ \si{kN/mm^2}$.
The critical energy release rate $G_c$ is defined as $G_c=2.7\ \si{N/mm}$. 
The loading increment is chosen as $\delta t=10^{-4}\ \si{s}$ and the
bulk regularization parameter $\kappa=10^{-10}$ is sufficiently small.
The mesh element diameter is set as $h_{\text{start}} = 0.088\ \si{mm}$. 
The end time $T$ is $0.0125\ \si{s}$, once the specimen is fully cracked.


\subsubsection{The single edge notched tension test}\label{config_tension}
The geometry and the material parameters are in line with the single edge
notched shear test in the previous section. The only difference is in the boundary conditions. As depicted in Figure~\ref{tension_geo}, it is pulled with a given force
in $y$-direction on the top boundary and the bottom boundary is fixed in $y$-direction. 
Also most of the numerical parameters are chosen as in the single edge notched shear test in Section~\ref{config_shear}. Here,
the mesh element diameter is set as $h_{\text{start}} = 0.044\ \si{mm}$.
In this setting, we use as loading increment $\delta t=10^{-5}\ \si{s}$, which
is necessary because of a very fast crack growth. The end time $T$ is
$0.00676\ \si{s}$, once the specimen is fully cracked.


\subsubsection{Quantities of Interest}
For both tests and all three studies, we discuss proper quantities of interest. In the first and second study we observe the load-displacement curves, where the load functions on the top boundary $\Gamma_{\text{top}}$ are computed by
\begin{align}
 \boldsymbol{\tau}=(F_x,F_y):=\int_{\Gamma_{\text{top}}} \boldsymbol{\sigma}(\boldsymbol{u}_h) \boldsymbol{n}\ \mathrm{d}s,\label{loading}
\end{align}
with the stress tensor $\boldsymbol{\sigma}(\boldsymbol{u}_h)$ and the normal vector $\boldsymbol{n}$. In the load-displacement curves the loading is displayed versus the displacements. Within the single edge notched shear test
we are particularly interested in the loading force $F_x$, within the tension test we are interested in the evaluation of $F_y$.

As a second quantity of interest, the bulk energy $E_b$ is defined as
\begin{align}
 E_b:=\int_{\Omega} \left( \left([1-\kappa]\varphi^2 + \kappa\right) \mu \operatorname{tr}\left(E_{\text{lin}}(\boldsymbol{u}_h)^2\right)+\frac{1}{2} \lambda \operatorname{tr}\left(E_{\text{lin}}(\boldsymbol{u}_h)\right)^2\right)\ \mathrm{d}x.\label{bulk}
\end{align}
Further the crack energy is computed via
\begin{align}
  E_c:=\frac{G_C}{2} \int_{\Omega}  \left(\frac{(\varphi-1)^2}{\epsilon} + \epsilon|\nabla \varphi|^2\right) \ \mathrm{d}x.\label{crack}
\end{align}
In addition, especially in Study $3$, we show snapshots of the phase-field function and the current adaptive mesh at certain time steps.


\subsection{Study $1$: uniform versus adaptive refinement}\label{study_one}
In a first study, the focus is on the comparison of adaptive meshes and a uniformly refined mesh. 
We compare the results of the load-displacement curves, the bulk and
the crack energy after one to six steps of adaptive refinement
starting with a coarser mesh than the used uniform refined
mesh for a fair comparison. 
First, the results of the single edge notched shear test are given.
For all tests executed in Study $1$, the relation between the discretization
parameter $h_{\text{start}}$ and $\epsilon$ is given by $\epsilon
=2h_{\text{start}}$ which means that we refine in $h$, but not in
$\epsilon$. The cell length $h_{\text{start}}$ is chosen as in the previous test.


\subsubsection{Results of the single edge notched shear test}

In the following, the load-displacement curves, the bulk and crack energy will be given for seven conducted numerical tests.
For the adaptive tests (named \textit{adaptive} $+$ number of refinement steps), the coarsest mesh is pre-refined three times, while the uniform computation is done on a mesh with six levels of global refinement.
In Figure~\ref{shear_study_one_loading}, seven load-displacement
curves are plotted. The loading $F_x$ is computed as defined
in~(\ref{loading}). In Figure~\ref{shear_study_one} the corresponding
maximal number of degrees of freedom, at each time step, are given for each
test. 
The tests called \textit{adaptive} $1$ to $6$ show the load-displacement curves computed on an adaptive mesh based
on a certain number of refinement cycles (one to six) according to Algorithm~\ref{alg:refinement}. The uniform mesh
consists of $66,820$ degrees of freedom ($6$ uniform refinement steps, $h=0.011\ \si{mm}$), 
which exceeds by far the number of degrees of freedom of the adaptive mesh after six steps of refinement based on the developed error estimator. 
\begin{figure}[htbp!]
\begin{minipage}{0.6\textwidth}
\begin{tikzpicture}[xscale=1,yscale=1]
\begin{axis}[
    ylabel = Load $F_x$ $\lbrack\si{N}\rbrack$,
    xlabel = Displacement $\lbrack\si{mm}\rbrack$,
 legend pos=north west, grid =major,
    x post scale = 1.5,
    y post scale = 1.5,
  xtick={0,0.002,0.004,0.006,0.008,0.01,0.012,0.014,0.016}, 
  ytick={0,100,200,300,400,500,600}
  ]
\addplot[red,densely dashed]
table[x=Loading,y=shear_2h_1,col sep=comma] {shear_loading_study_one.csv};
\addlegendentry{adaptive 1}
\addplot[green]
table[x=Loading,y=shear_2h_2,col sep=comma] {shear_loading_study_one.csv};
\addlegendentry{adaptive 2}
\addplot[violet,densely dotted]
table[x=Loading,y=shear_2h_3,col sep=comma] {shear_loading_study_one.csv};
\addlegendentry{adaptive 3}
\addplot[orange]
table[x=Loading,y=shear_2h_4,col sep=comma] {shear_loading_study_one.csv};
\addlegendentry{adaptive 4}
\addplot[magenta,dash dot dot]
table[x=Loading,y=shear_2h_5,col sep=comma] {shear_loading_study_one.csv};
\addlegendentry{adaptive 5}
\addplot[cyan]
table[x=Loading,y=shear_2h_6,col sep=comma] {shear_loading_study_one.csv};
\addlegendentry{adaptive 6}
\addplot[blue,dashed]
table[x=Loading,y=shear_2h_uniform,col sep=comma] {shear_loading_uniform.csv}; 
\addlegendentry{uniform}
\end{axis}
\end{tikzpicture}
\caption{Load-displacement curves for the single edge notched shear test with six steps of uniform refinement compared to six steps of adaptive refinement and three steps of pre-refinement. 
We fix $\epsilon=2h_{\text{start}}$ with $h_{\text{start}}$ as the discretization length on the coarsest mesh. It means, $\epsilon$ is fixed while $h$ decreases over the adaptive refinement steps.}\label{shear_study_one_loading}
\end{minipage}
\hfill
\begin{minipage}{0.23\textwidth}
\renewcommand*{\arraystretch}{1.2}
\begin{tabular}{|c|r|}\hline
\multicolumn{1}{|c}{Test name}  & \multicolumn{1}{|c|}{$\#$ DoFs} \\ \hline \hline
   adaptive $1$  & $2,672$\\ \hline 
   adaptive $2$  & $3,928$ \\ \hline 
   adaptive $3$  & $9,836$ \\ \hline 
   adaptive $4$  & $12,740$ \\ \hline 
    adaptive $5$  &  $16,024$ \\ \hline
   adaptive $6$  & $18,196$ \\ \hline 
   uniform &  $66,820$ \\ \hline 
 \end{tabular}
\caption{Maximal number of degrees of freedom for each test case for Study $1$ designed for the single edge notched shear test.}\label{shear_study_one}
\end{minipage}
\end{figure}

In Figure~\ref{shear_study_one_loading}, with 
only $18,196$ degrees of freedom (test \textit{adaptive} $6$) we obtain nearly the same load-displacement curve as on the uniform refined mesh with $66,820$ degrees of freedom. 
Furthermore, with an increasing number of refinement steps, we observe convergence
towards the load-displacement curve computed with uniform mesh refinement.

We continue with the quantities bulk energy (defined in Equation~(\ref{bulk})) and crack energy (defined in Equation~(\ref{crack})) with the same test setup as 
listed in Figure~\ref{shear_study_one}. Notice that neither the
load in Figure~\ref{shear_study_one_loading} nor the
bulk-energy~\ref{shear_study_one_bulk} vanish at the end of the
loading process. This is due to the fact, cf.~\cite{AmGeraLoren15},
that the chosen stress-splitting does not allow for complete fracture. 

\begin{figure}[htbp!]
\begin{minipage}{0.48\textwidth}
\begin{tikzpicture}[xscale=0.75,yscale=0.75]
\begin{axis}[
    ylabel = Bulk energy $E_b$ $\lbrack\si{J}\rbrack$,
    xlabel = Displacement $\lbrack\si{mm}\rbrack$,
 legend pos=north west, grid =major,
    x post scale = 1.4,
    y post scale = 1.5,
  xtick={0,0.005,0.01,0.015,0.02}, 
  ytick={0,0.25,0.5,0.75,1.0,1.25,1.5,1.75,2.0}
  ]
\addplot[red,densely dashed]
table[x=Loading,y=shear_2h_1,col sep=comma] {shear_bulk_study_one.csv};
\addlegendentry{adaptive 1}
\addplot[green]
table[x=Loading,y=shear_2h_2,col sep=comma] {shear_bulk_study_one.csv};
\addlegendentry{adaptive 2}
\addplot[violet,densely dotted]
table[x=Loading,y=shear_2h_3,col sep=comma] {shear_bulk_study_one.csv};
\addlegendentry{adaptive 3}
\addplot[orange]
table[x=Loading,y=shear_2h_4,col sep=comma] {shear_bulk_study_one.csv};
\addlegendentry{adaptive 4}
\addplot[magenta,dash dot dot]
table[x=Loading,y=shear_2h_5,col sep=comma] {shear_bulk_study_one.csv};
\addlegendentry{adaptive 5}
\addplot[cyan]
table[x=Loading,y=shear_2h_6,col sep=comma] {shear_bulk_study_one.csv};
\addlegendentry{adaptive 6}
\addplot[blue,dashed]
table[x=Loading,y=shear_2h_bulk_uniform,col sep=comma] {shear_bulk_uniform.csv};
\addlegendentry{uniform}
\end{axis}
\end{tikzpicture}
\caption{Bulk energy for the single edge notched shear test with six steps of uniform refinement compared to six steps of adaptive refinement and three steps of pre-refinement. 
We fix $\epsilon=2h_{\text{start}}$, while $h$ decreases over the adaptive refinement steps.}\label{shear_study_one_bulk}
\end{minipage}
\hfill
\begin{minipage}{0.48\textwidth}
\begin{tikzpicture}[xscale=0.75,yscale=0.75]
\begin{axis}[
    ylabel = Crack energy $E_c$ $\lbrack\si{J}\rbrack$,
    xlabel = Displacement $\lbrack\si{mm}\rbrack$,
 legend pos=north west, grid =major,
    x post scale = 1.4,
    y post scale = 1.5,
  xtick={0,0.005,0.01,0.015,0.02}, 
  ytick={0,0.25,0.5,0.75,1.0,1.25,1.5,1.75,2.0}
  ]
\addplot[red,densely dashed]
table[x=Loading,y=shear_2h_1,col sep=comma] {shear_crack_study_one.csv};
\addlegendentry{adaptive 1}
\addplot[green]
table[x=Loading,y=shear_2h_2,col sep=comma] {shear_crack_study_one.csv};
\addlegendentry{adaptive 2}
\addplot[violet,densely dotted]
table[x=Loading,y=shear_2h_3,col sep=comma] {shear_crack_study_one.csv};
\addlegendentry{adaptive 3}
\addplot[orange]
table[x=Loading,y=shear_2h_4,col sep=comma] {shear_crack_study_one.csv};
\addlegendentry{adaptive 4}
\addplot[magenta,dash dot dot]
table[x=Loading,y=shear_2h_5,col sep=comma] {shear_crack_study_one.csv};
\addlegendentry{adaptive 5}
\addplot[cyan]
table[x=Loading,y=shear_2h_6,col sep=comma] {shear_crack_study_one.csv};
\addlegendentry{adaptive 6}
\addplot[blue,dashed]
table[x=Loading,y=shear_2h_crack_uniform,col sep=comma] {shear_crack_uniform.csv};
\addlegendentry{uniform}
\end{axis}
\end{tikzpicture}
\caption{Crack energy for the single edge notched shear test with six steps of uniform refinement compared to six steps of adaptive refinement and three steps of pre-refinement. 
We fix $\epsilon=2h_{\text{start}}$, while $h$ decreases over the adaptive refinement steps.}\label{shear_study_one_crack}
\end{minipage}
\end{figure}

Comparing the results of the bulk and crack energy depicted in the Figures~\ref{shear_study_one_bulk} and~\ref{shear_study_one_crack}, we observe that, once crack propagation starts, 
the bulk energy decreases and the crack energy increases.
Secondly, similar to the load-displacement curves the course of the bulk and crack energy over time/loading 
tends towards the curves given by the computation on a uniformly refined mesh with an increasing number of refinement steps.


\subsubsection{Results of the single edge notched tension test}
In the following, the load-displacement curves, the bulk and the crack energy are computed and presented for six conducted numerical tests based on the tension test.
For this test, the coarsest mesh is pre-refined four times, while the uniform
computation is done as in the shear test on a mesh with six levels of global refinement.

In Figure~\ref{tension_study_one_loading}, six load-displacement curves are plotted. In Figure~\ref{shear_study_one}, the corresponding maximal number of degrees of freedom are given for each test. 
As in the shear example, the tests called \textit{adaptive} $1$ to $5$ show the load-displacement curves (loading $F_y$ computed via~(\ref{loading})) based on an adaptive mesh with a certain number of refinement steps ($1$ to $5$). The uniform mesh consists of
$66,820$ degrees of freedom ($6$ uniform refinement steps, $h=0.011\ \si{mm}$), which is more than the adaptive mesh after five steps of refinement based on the developed error estimator.
Comparing the test based on a uniform mesh and the one with three adaptive refinement steps 
(dotted violet curve in Figure ~\ref{tension_study_one_loading}), the load-displacement curves are quite similar, although the adaptive test computes with maximal $18,408$ degrees of freedom per time step.

\begin{figure}[htbp!]
\begin{minipage}{0.6\textwidth}
\centering
\begin{tikzpicture}[xscale=1,yscale=1]
\begin{axis}[
    ylabel = Load $F_x$ $\lbrack\si{N}\rbrack$,
    xlabel = Displacement $\lbrack\si{mm}\rbrack$,
 legend pos=north west, grid =major,
    x post scale = 1.5,
    y post scale = 1.5,
  xtick={0,0.001,0.002,0.003,0.004,0.005,0.006,0.007,0.008}, 
  ytick={0,100,200,300,400,500,600,700,800}
  ]
 \addplot[red,densely dashed]
 table[x=Loading,y=tension_2h_1,col sep=comma] {tension_loading_study_one.csv};
 \addlegendentry{adaptive 1}
 \addplot[green]
 table[x=Loading,y=tension_2h_2,col sep=comma] {tension_loading_study_one.csv};
 \addlegendentry{adaptive 2}
 \addplot[violet,densely dotted]
 table[x=Loading,y=tension_2h_3,col sep=comma] {tension_loading_study_one.csv};
 \addlegendentry{adaptive 3}
 \addplot[orange]
 table[x=Loading,y=tension_2h_4,col sep=comma] {tension_loading_study_one.csv};
 \addlegendentry{adaptive 4}
 \addplot[magenta,dash dot dot]
 table[x=Loading,y=tension_2h_5,col sep=comma] {tension_loading_study_one.csv};
 \addlegendentry{adaptive 5}
\addplot[blue,dashed]
table[x=Loading,y=tension_2h_uniform,col sep=comma] {tension_loading_uniform.csv};
\addlegendentry{uniform}
\end{axis}
\end{tikzpicture}
\caption{Load-displacement curves for the single edge notched tension test with six steps of uniform refinement compared to five steps of adaptive refinement and four steps of pre-refinement. 
We fix $\epsilon=2h_{\text{start}}$, while $h$ decreases over the adaptive refinement steps.}\label{tension_study_one_loading}
\end{minipage}
\hfill
\begin{minipage}{0.23\textwidth}
\renewcommand*{\arraystretch}{1.2}
\begin{tabular}{|c|r|}\hline
\multicolumn{1}{|c}{Test name}  & \multicolumn{1}{|c|}{$\#$ DoFs} \\ \hline \hline
   adaptive 1  & $6,468$\\ \hline
   adaptive 2  & $11,052$ \\ \hline 
   adaptive 3  & $18,408$ \\ \hline 
   adaptive 4  & $37,916$ \\ \hline 
    adaptive 5  &  $61,704$ \\ \hline 
   uniform &  $66,820$ \\ \hline
 \end{tabular}
\caption{Maximal number of degrees of freedom for each test case for Study $1$ designed for the single edge notched tension test.}\label{tension_study_one}
\end{minipage}
\end{figure}

\begin{figure}[htbp!]
\begin{minipage}{0.48\textwidth}
\begin{tikzpicture}[xscale=0.75,yscale=0.75]
\begin{axis}[
    ylabel = Bulk energy $E_b$ $\lbrack\si{J}\rbrack$,
    xlabel = Displacement $\lbrack\si{mm}\rbrack$,
 legend pos=north west, grid =major,
    x post scale = 1.4,
    y post scale = 1.5,
   xtick={0,0.001,0.002,0.003,0.004,0.005,0.006,0.007,0.008,0.009}, 
  ytick={0,0.25,0.5,0.75,1.0,1.25,1.5}
  ]
\addplot[red,densely dashed]
table[x=Loading,y=tension_2h_1,col sep=comma] {tension_bulk_study_one.csv};
\addlegendentry{adaptive 1}
\addplot[green]
table[x=Loading,y=tension_2h_2,col sep=comma] {tension_bulk_study_one.csv};
\addlegendentry{adaptive 2}
\addplot[violet, densely dotted]
table[x=Loading,y=tension_2h_3,col sep=comma] {tension_bulk_study_one.csv};
\addlegendentry{adaptive 3}
\addplot[orange]
table[x=Loading,y=tension_2h_4,col sep=comma] {tension_bulk_study_one.csv};
\addlegendentry{adaptive 4}
\addplot[magenta,dash dot dot]
table[x=Loading,y=tension_2h_5,col sep=comma] {tension_bulk_study_one.csv};
\addlegendentry{adaptive 5}
\addplot[blue,dashed]
table[x=Loading,y=tension_2h_bulk_uniform,col sep=comma] {tension_bulk_uniform.csv};
\addlegendentry{uniform}
\end{axis}
\end{tikzpicture}
\caption{Bulk energy for the single edge notched tension test with six steps of uniform refinement compared to five steps of adaptive refinement and four steps of pre-refinement. 
We fix $\epsilon=2h_{\text{start}}$, while $h$ decreases over the adaptive refinement steps.}\label{tension_study_one_bulk}
\end{minipage}
\hfill
\begin{minipage}{0.48\textwidth}
\begin{tikzpicture}[xscale=0.75,yscale=0.75]
\begin{axis}[
    ylabel = Crack energy $E_c$ $\lbrack\si{J}\rbrack$,
    xlabel = Displacement $\lbrack\si{mm}\rbrack$,
 legend pos=north west, grid =major,
    x post scale = 1.4,
    y post scale = 1.5,
  xtick={0,0.001,0.002,0.003,0.004,0.005,0.006,0.007,0.008,0.009}, 
  ytick={0,0.25,0.5,0.75,1.0,1.25,1.5}
  ]
\addplot[red,densely dashed]
table[x=Loading,y=tension_2h_1,col sep=comma] {tension_crack_study_one.csv};
\addlegendentry{adaptive 1}
\addplot[green]
table[x=Loading,y=tension_2h_2,col sep=comma] {tension_crack_study_one.csv};
\addlegendentry{adaptive 2}
\addplot[violet,densely dotted]
table[x=Loading,y=tension_2h_3,col sep=comma] {tension_crack_study_one.csv};
\addlegendentry{adaptive 3}
\addplot[orange]
table[x=Loading,y=tension_2h_4,col sep=comma] {tension_crack_study_one.csv};
\addlegendentry{adaptive 4}
\addplot[magenta,dash dot dot]
table[x=Loading,y=tension_2h_5,col sep=comma] {tension_crack_study_one.csv};
\addlegendentry{adaptive 5}
\addplot[blue,dashed]
table[x=Loading,y=tension_2h_crack_uniform,col sep=comma] {tension_crack_uniform.csv};
\addlegendentry{uniform}
\end{axis}
\end{tikzpicture}
\caption{Crack energy for the single edge notched tension test with six steps of uniform refinement compared to five steps of adaptive refinement and four steps of pre-refinement. 
We fix $\epsilon=2h_{\text{start}}$, while $h$ decreases over the adaptive refinement steps.}\label{tension_study_one_crack}
\end{minipage}
\end{figure}
  
As already discussed for the single edge notched shear test, the evolution of the bulk and crack energy visualized in the Figures~\ref{tension_study_one_bulk} and~\ref{tension_study_one_crack} is of the expected behavior: as long as the crack does not grow, only the bulk energy increases. 
After crack propagation, bulk energy is dissipated into crack energy.


\subsection{Study 2: $\epsilon-h$ relationship}\label{study_two}
Now, we are interested in the numerical results of tests with
different relations of the discretization parameter $h_{\text{start}}$ on the initial mesh and the crack bandwidth $\epsilon$.
We investigate the $h_{\text{start}}-\epsilon$ relationship as follows:
\begin{itemize}
 \item Case 1: $\epsilon =2 h_{\text{start}}$ with a fixed discretization parameter $h_{\text{start}} = 0.088\ \si{mm}$, while $h$ is decreasing during increasing refinement steps
 \item  Case 2: $\epsilon = h_{\text{start}}$ with $h_{\text{start}} = 0.088\ \si{mm}$ fixed during refinement
 \item Case 3: $\epsilon =\frac{1}{2} h_{\text{start}}$ with $h_{\text{start}} = 0.088\ \si{mm}$ fixed during refinement
\end{itemize}


\subsubsection{Results of the single edge notched shear test}
\begin{figure}[htbp!]
\centering
\begin{tikzpicture}[xscale=0.75,yscale=0.75]
\begin{axis}[
    ylabel = Load $F_x$ $\lbrack\si{N}\rbrack$,
    xlabel = Displacement $\lbrack\si{mm}\rbrack$,
 legend pos=north west, grid =major,
    x post scale = 1.75,
    y post scale = 1.5,
  xtick={0,0.005,0.01,0.015,0.02}, 
  ytick={0,100,200,300,400,500,600}
  ]
\addplot[blue,dashed]  
table[x=Loading,y=shear_2h,col sep=comma] {shear_loading.csv};
\addlegendentry{$\epsilon=2h_{\text{start}}$ adaptive}
\addplot[violet,densely dotted]  
table[x=Loading,y=shear_h,col sep=comma] {shear_loading.csv};
\addlegendentry{$\epsilon=h_{\text{start}}$ adaptive}
\addplot[orange]  
table[x=Loading,y=shear_05h,col sep=comma] {shear_loading.csv};
\addlegendentry{$\epsilon=0.5h_{\text{start}}$ adaptive}
\end{axis}
\end{tikzpicture}
\caption{Load-displacement curves for the single edge notched shear test with three steps of global pre-refinement and six steps of adaptive refinement and $\epsilon=2h_{\text{start}}$, $\epsilon=h_{\text{start}}$ and $\epsilon=0.5 h_{\text{start}}$.}\label{shear_adaptive_eps_load}
\end{figure}
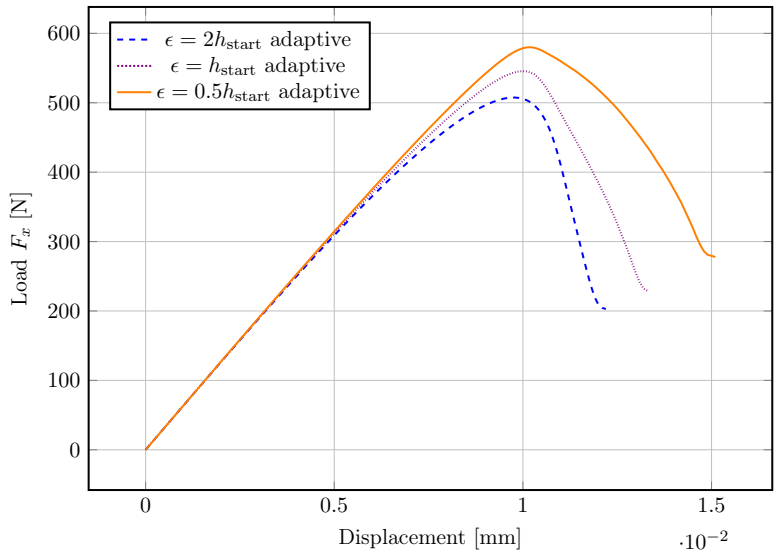

In Figure~\ref{shear_adaptive_eps_load}, the load-displacement curves are depicted for the three test cases mentioned above. 
For all load-displacement curves one can observe convergence, even if for $\epsilon \geq h_{\text{start}}$, the results are more precise and more similar to the one in Figure~\ref{shear_study_one_loading}.

The same can be observed for the bulk and crack energy in the Figures~\ref{shear_adaptive_eps_bulk} and~\ref{shear_adaptive_eps_crack}. 
The general course of the energies is acceptable, but with a decreasing
bandwidth $\epsilon$, the curves are less accurate against the starting point
of crack propagation. Similar observations were made in~\cite{HeWheWi15}.

\begin{figure}[htbp!]
\begin{minipage}{0.48\textwidth}
\begin{tikzpicture}[xscale=0.75,yscale=0.75]
\begin{axis}[
    ylabel = Bulk energy $E_b$ $\lbrack\si{J}\rbrack$,
    xlabel = Displacement $\lbrack\si{mm}\rbrack$,
 legend pos=north west, grid =major,
    x post scale = 1.4,
    y post scale = 1.5,
  xtick={0,0.005,0.01,0.015,0.02}, 
  ytick={0,0.5,1.0,1.5,2.0}
  ]
\addplot[blue,dashed]  
table[x=Loading,y=shear_2h_bulk,col sep=comma] {shear_bulk.csv};
\addlegendentry{$\epsilon=2h$ adaptive}
\addplot[violet,densely dotted]  
table[x=Loading,y=shear_h_bulk,col sep=comma] {shear_bulk.csv};
\addlegendentry{$\epsilon=h$ adaptive}
\addplot[orange]  
table[x=Loading,y=shear_05h_bulk,col sep=comma] {shear_bulk.csv};
\addlegendentry{$\epsilon=0.5h$ adaptive}
\end{axis}
\end{tikzpicture}
\caption{Bulk energy for the single edge notched shear test with three steps of global pre-refinement and six steps of adaptive refinement and $\epsilon=2h_{\text{start}}$, $\epsilon=h_{\text{start}}$ and $\epsilon=0.5 h_{\text{start}}$.}\label{shear_adaptive_eps_bulk}
\end{minipage}
\hfill
\begin{minipage}{0.48\textwidth}
\begin{tikzpicture}[xscale=0.75,yscale=0.75]
\begin{axis}[
    ylabel = Crack energy $E_c$ $\lbrack\si{J}\rbrack$,
    xlabel = Displacement $\lbrack\si{mm}\rbrack$,
 legend pos=north west, grid =major,
    x post scale = 1.4,
    y post scale = 1.5,
  xtick={0,0.005,0.01,0.015,0.02}, 
  ytick={0,0.5,1.0,1.5,2.0}
  ]
\addplot[blue,dashed]  
table[x=Loading,y=shear_2h_crack,col sep=comma] {shear_crack.csv};
\addlegendentry{$\epsilon=2h_{\text{start}}$ adaptive}
\addplot[violet,densely dotted]  
table[x=Loading,y=shear_h_crack,col sep=comma] {shear_crack.csv};
\addlegendentry{$\epsilon=h_{\text{start}}$ adaptive}
\addplot[orange]  
table[x=Loading,y=shear_05h_crack,col sep=comma] {shear_crack.csv};
\addlegendentry{$\epsilon=0.5h_{\text{start}}$ adaptive}
\end{axis}
\end{tikzpicture}
\caption{Crack energy for the single edge notched shear test with three steps of global pre-refinement and six steps of adaptive refinement and $\epsilon=2h_{\text{start}}$, $\epsilon=h_{\text{start}}$ and $\epsilon=0.5 h_{\text{start}}$.}\label{shear_adaptive_eps_crack}
\end{minipage}
\end{figure}


\subsubsection{Results of the single edge notched tension test}
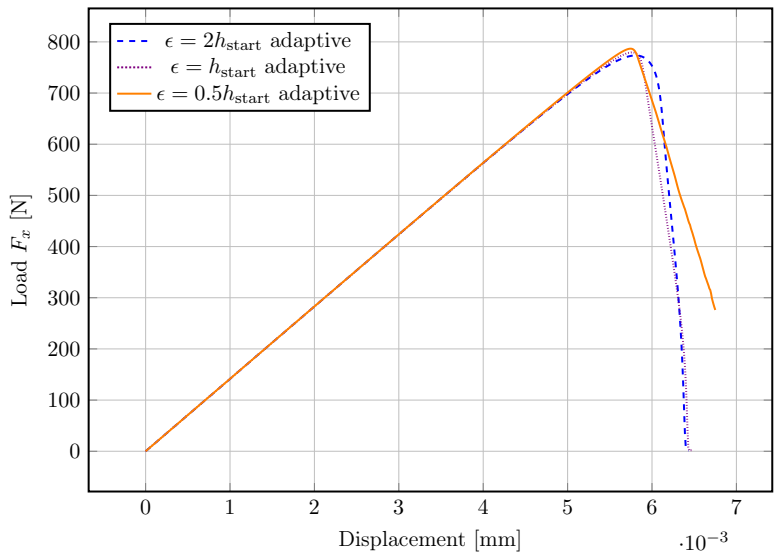
\begin{figure}[htbp!]
\centering
\begin{tikzpicture}[xscale=0.75,yscale=0.75]
\begin{axis}[
    ylabel = Load $F_x$ $\lbrack\si{N}\rbrack$,
    xlabel = Displacement $\lbrack\si{mm}\rbrack$,
 legend pos=north west, grid =major,
    x post scale = 1.75,
    y post scale = 1.5,
  xtick={0,0.001,0.002,0.003,0.004,0.005,0.006,0.007,0.008,0.009}, 
  ytick={0,100,200,300,400,500,600,700,800,900}
  ]
\addplot[blue,dashed]  
table[x=Loading,y=tension_2h,col sep=comma] {tension_loading.csv};
\addlegendentry{$\epsilon=2h_{\text{start}}$ adaptive}
\addplot[violet,densely dotted]  
table[x=Loading,y=tension_h,col sep=comma] {tension_loading.csv};
\addlegendentry{$\epsilon=h_{\text{start}}$ adaptive}
\addplot[orange]  
table[x=Loading,y=tension_05h,col sep=comma] {tension_loading.csv};
\addlegendentry{$\epsilon=0.5h_{\text{start}}$ adaptive}
\end{axis}
\end{tikzpicture}
\caption{Load-displacement curves for the single edge notched tension test with four steps of global pre-refinement and five steps of adaptive refinement and $\epsilon=2h_{\text{start}}$, $\epsilon=h_{\text{start}}$ and $\epsilon=0.5 h_{\text{start}}$.}\label{tension_adaptive_eps_load}
\end{figure}

As expected, the load-displacement curves in Figure~\ref{tension_adaptive_eps_load} for the three test cases $\epsilon=2h_{\text{start}}$, $\epsilon=h_{\text{start}}$  
and $\epsilon=0.5h_{\text{start}}$ show an increase of the loading with an increasing
displacement and a steep descent at the point where the crack starts
propagating. For $\epsilon=0.5h_{\text{start}}$ apparently the transition zone can be resolved sufficiently after five steps of adaptive refinement
based on the initial mesh. This behavior can be recognized in the plotted energies in the Figures~\ref{tension_adaptive_eps_bulk} and~\ref{tension_adaptive_eps_crack}.

\begin{figure}[htbp!]
\begin{minipage}{0.48\textwidth}
\begin{tikzpicture}[xscale=0.75,yscale=0.75]
\begin{axis}[
    ylabel = Bulk energy $E_b$ $\lbrack\si{J}\rbrack$,
    xlabel = Displacement $\lbrack\si{mm}\rbrack$,
 legend pos=north west, grid =major,
    x post scale = 1.4,
    y post scale = 1.5,
  xtick={0,0.001,0.002,0.003,0.004,0.005,0.006,0.007,0.008,0.009}, 
  ytick={0,0.5,1.0,1.5,2.0}
  ]
\addplot[blue,dashed]  
table[x=Loading,y=tension_2h_bulk,col sep=comma] {tension_bulk.csv};
\addlegendentry{$\epsilon=2h_{\text{start}}$ adaptive}
\addplot[violet,densely dotted]  
table[x=Loading,y=tension_h_bulk,col sep=comma] {tension_bulk.csv};
\addlegendentry{$\epsilon=h_{\text{start}}$ adaptive}
\addplot[orange]  
table[x=Loading,y=tension_05h_bulk,col sep=comma] {tension_bulk.csv};
\addlegendentry{$\epsilon=0.5h_{\text{start}}$ adaptive}
\end{axis}
\end{tikzpicture}
\caption{Bulk energy for the single edge notched tension test with four steps of global pre-refinement and five steps of adaptive refinement and $\epsilon=2h_{\text{start}}$, $\epsilon=h_{\text{start}}$ and $\epsilon=0.5 h_{\text{start}}$.}\label{tension_adaptive_eps_bulk}
\end{minipage}
\hfill
\begin{minipage}{0.48\textwidth}
\begin{tikzpicture}[xscale=0.75,yscale=0.75]
\begin{axis}[
    ylabel = Crack energy $E_c$ $\lbrack\si{J}\rbrack$,
    xlabel = Displacement $\lbrack\si{mm}\rbrack$,
 legend pos=north west, grid =major,
    x post scale = 1.4,
    y post scale = 1.5,
  xtick={0,0.001,0.002,0.003,0.004,0.005,0.006,0.007,0.008,0.009}, 
  ytick={0,0.5,1.0,1.5,2.0}
  ]
\addplot[blue,dashed]  
table[x=Loading,y=tension_2h_crack,col sep=comma] {tension_crack.csv};
\addlegendentry{$\epsilon=2h_{\text{start}}$ adaptive}
\addplot[violet,densely dotted]  
table[x=Loading,y=tension_h_crack,col sep=comma] {tension_crack.csv};
\addlegendentry{$\epsilon=h_{\text{start}}$ adaptive}
\addplot[orange]  
table[x=Loading,y=tension_05h_crack,col sep=comma] {tension_crack.csv};
\addlegendentry{$\epsilon=0.5h_{\text{start}}$ adaptive}
\end{axis}
\end{tikzpicture}
\caption{Crack energy for the single edge notched tension test with four steps of global pre-refinement and five steps of adaptive refinement and $\epsilon=2h_{\text{start}}$, $\epsilon=h_{\text{start}}$ and $\epsilon=0.5 h_{\text{start}}$.}\label{tension_adaptive_eps_crack}
\end{minipage}
\end{figure}

The results in Section~\ref{study_two}, in particular
$\epsilon=0.5h_{\text{start}}$, indicate that our adaptive algorithm
performs well even if the relation $\epsilon > 2h_{\text{start}}$ is
not satisfied as it is required on uniform meshes for a reasonable
accuracy of the crack approximation, cf.~\cite[Figure~4]{MieWelHof10a}.


\subsection{Study 3: adaptive meshes and local error indicators}\label{study_three}
In this final study, we present snapshots of the adaptively refined meshes 
and provide visualizations of the local error indicators.


\subsubsection{Results of the single edge notched shear test}

\begin{figure}[htbp!]
\begin{minipage}{0.3\textwidth}
 \includegraphics[width=0.95\textwidth]{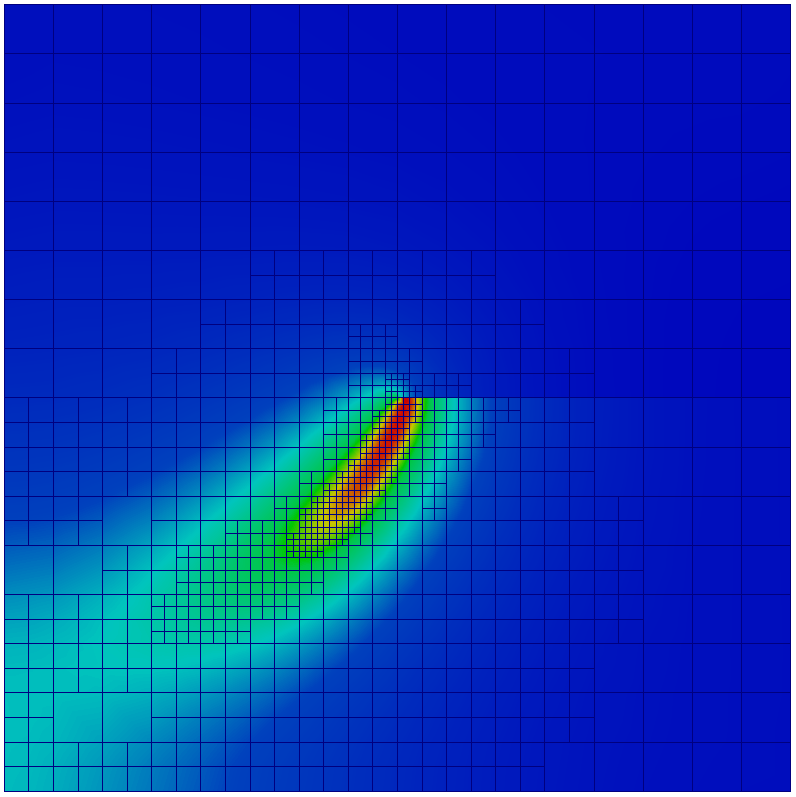}
\end{minipage}
\hfill
\begin{minipage}{0.3\textwidth}
 \includegraphics[width=0.95\textwidth]{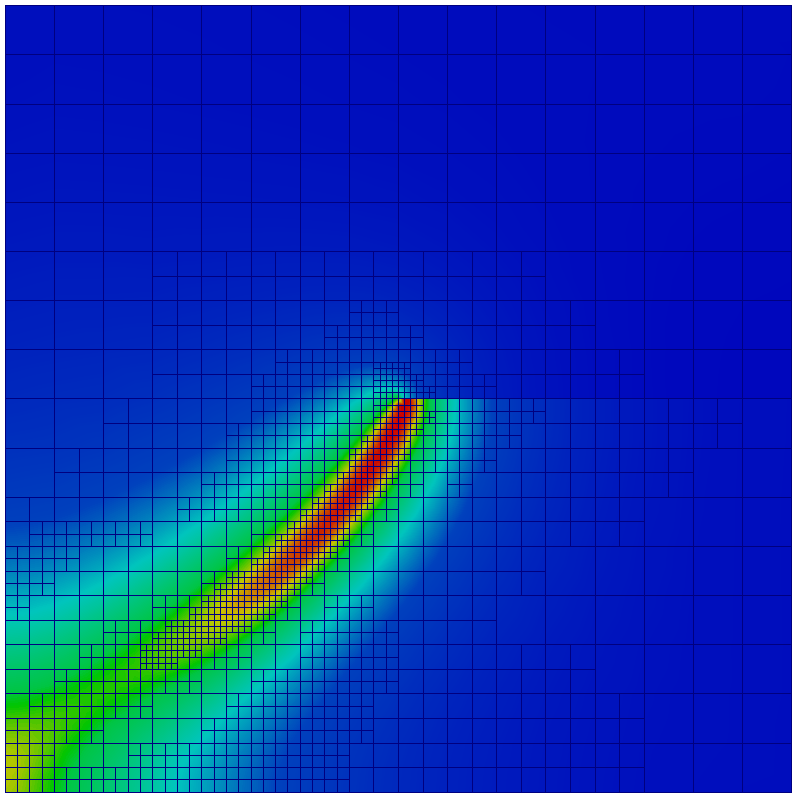}
 \end{minipage}
 \hfill
\begin{minipage}{0.375\textwidth}
 \includegraphics[width=0.95\textwidth]{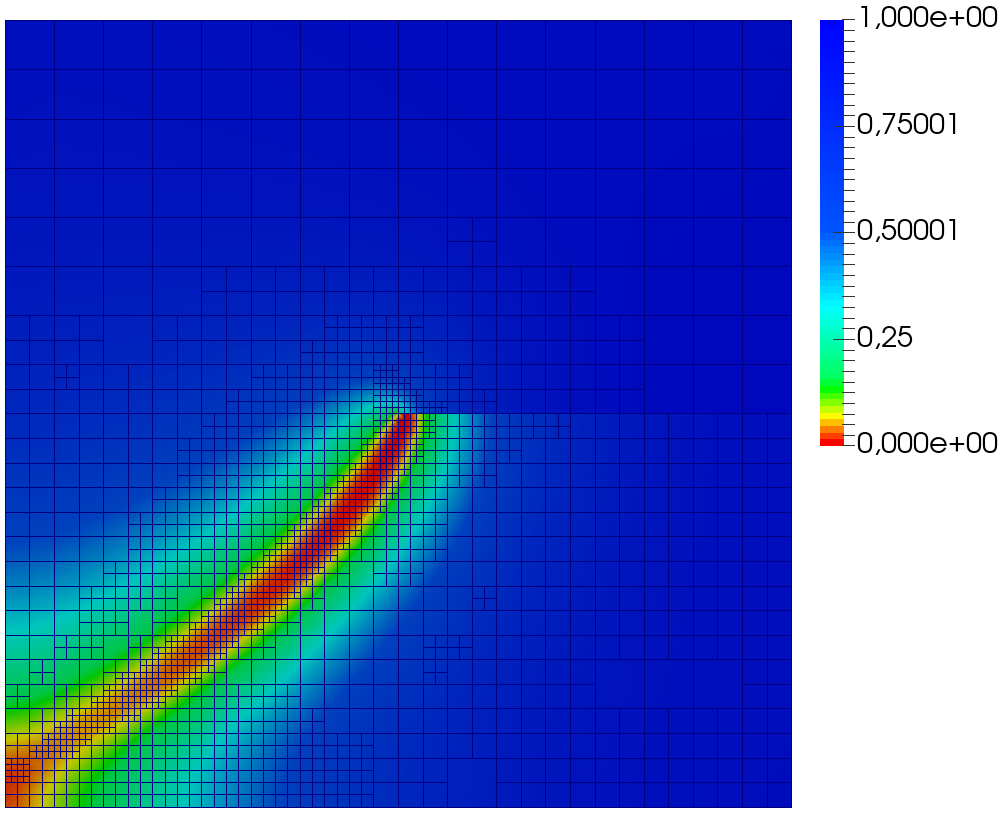}
 \end{minipage}\\
\vfill
\begin{minipage}{0.3\textwidth}
 \includegraphics[width=0.95\textwidth]{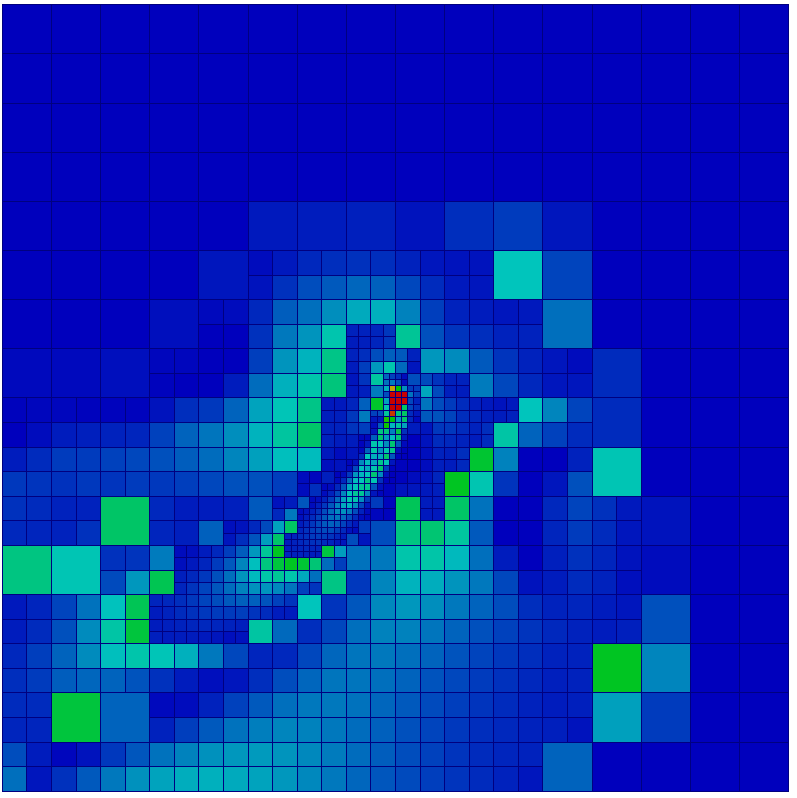}
\end{minipage}
\hfill
\begin{minipage}{0.3\textwidth}
 \includegraphics[width=0.95\textwidth]{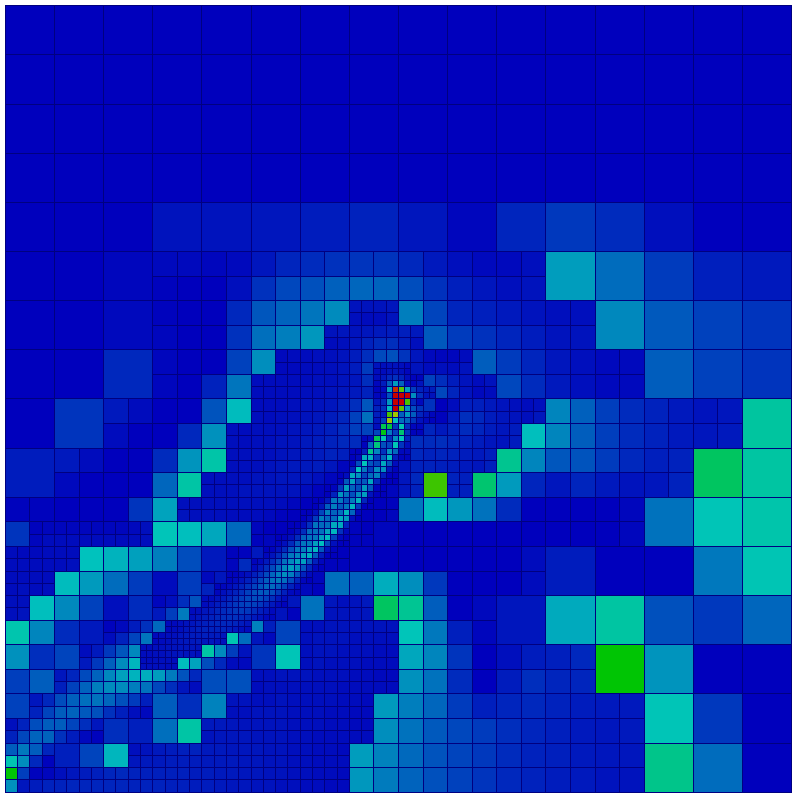}
 \end{minipage}
 \hfill
\begin{minipage}{0.375\textwidth}
 \includegraphics[width=0.95\textwidth]{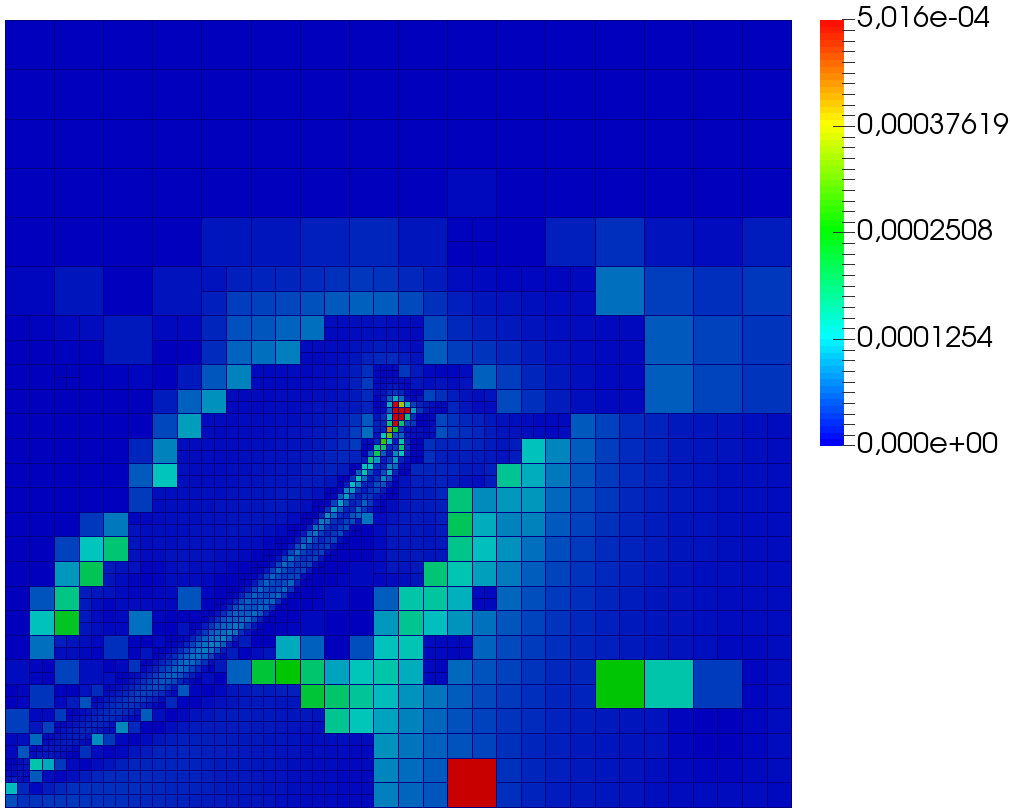}
 \end{minipage}
\caption{The phase-field function and the error indicators, respectively, after three refinement steps given in certain time steps (after $0.0115$, $0.00119$ and $0.00125\ \si{s}$) for the single edge notched shear test given on the current adaptive mesh to visualize the refinement strategy.}\label{phasefield_shear_3}
\end{figure}

\begin{figure}[htbp!]
\begin{minipage}{0.3\textwidth}
 \includegraphics[width=0.95\textwidth]{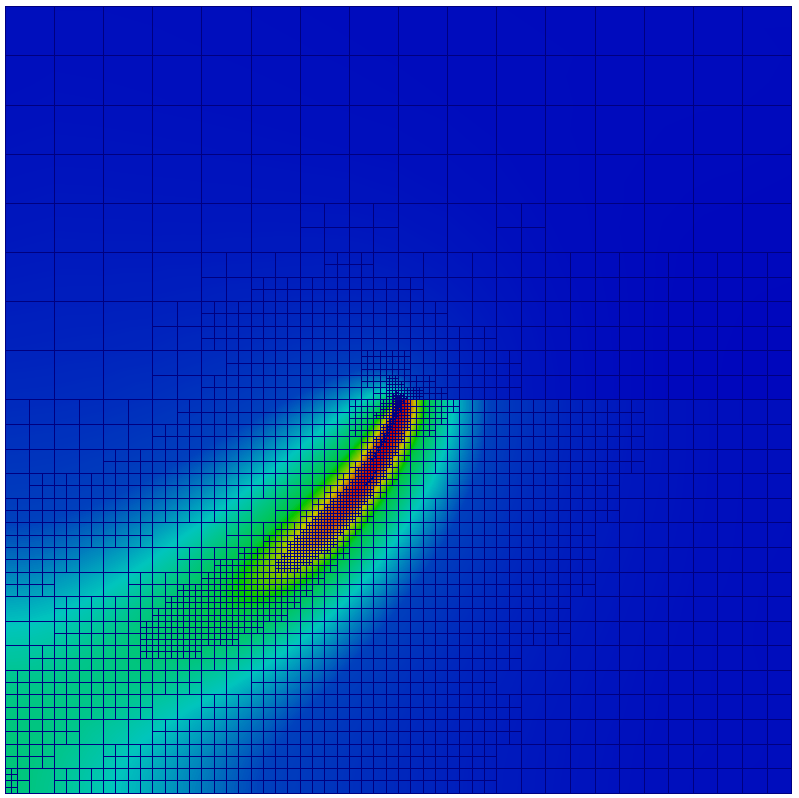}
\end{minipage}
\hfill
\begin{minipage}{0.3\textwidth}
 \includegraphics[width=0.95\textwidth]{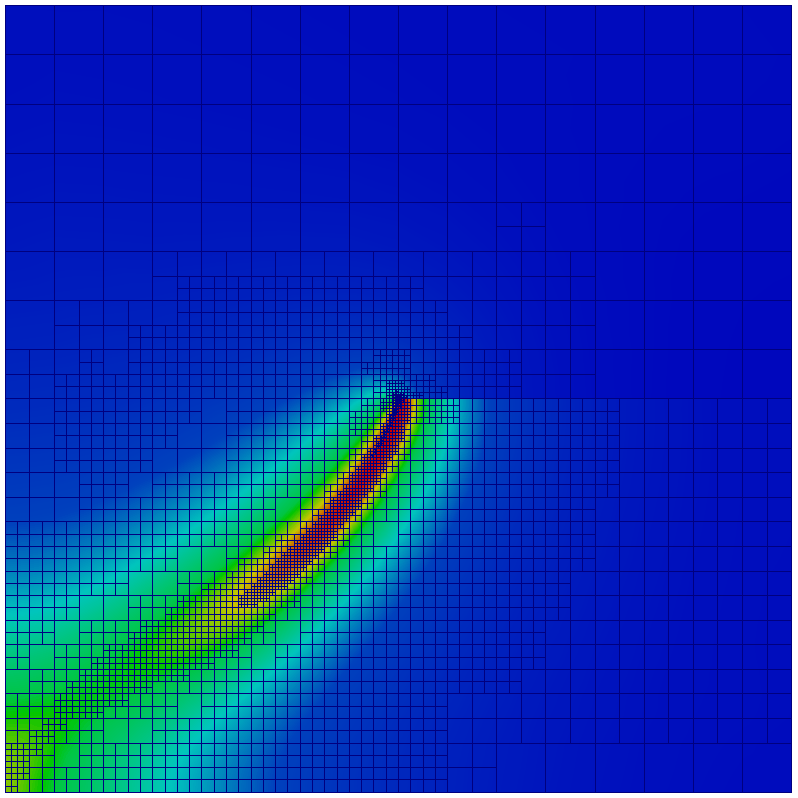}
 \end{minipage}
 \hfill
\begin{minipage}{0.375\textwidth}
 \includegraphics[width=0.95\textwidth]{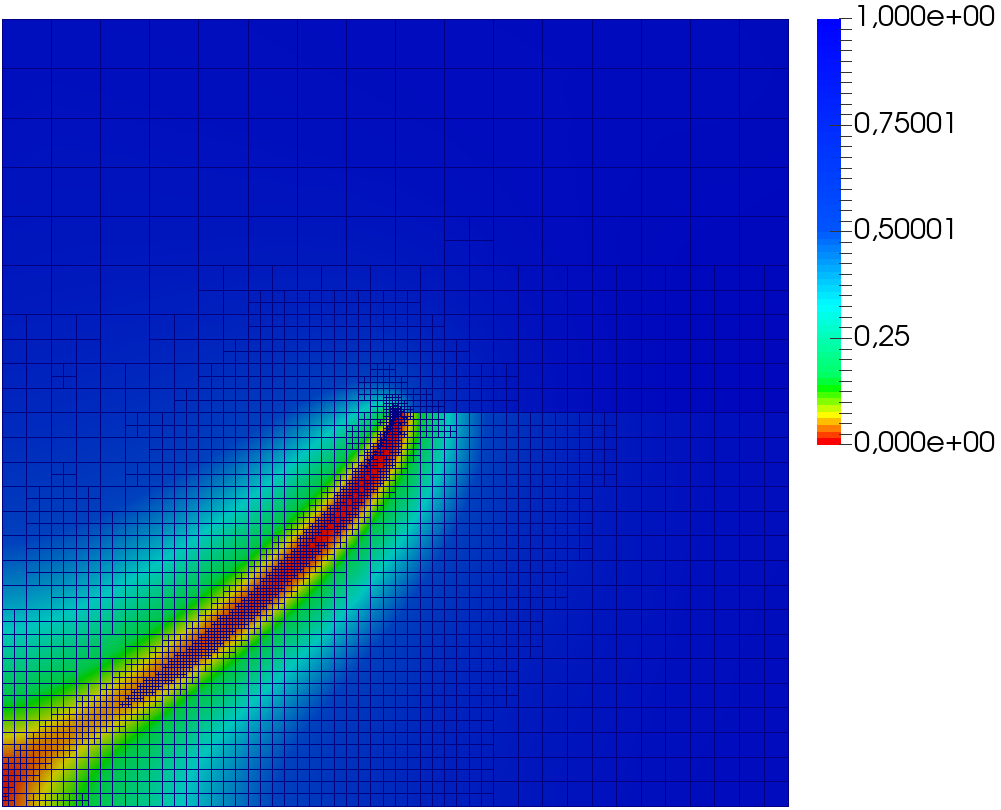}
 \end{minipage}\\
 \vfill
\begin{minipage}{0.3\textwidth}
 \includegraphics[width=0.95\textwidth]{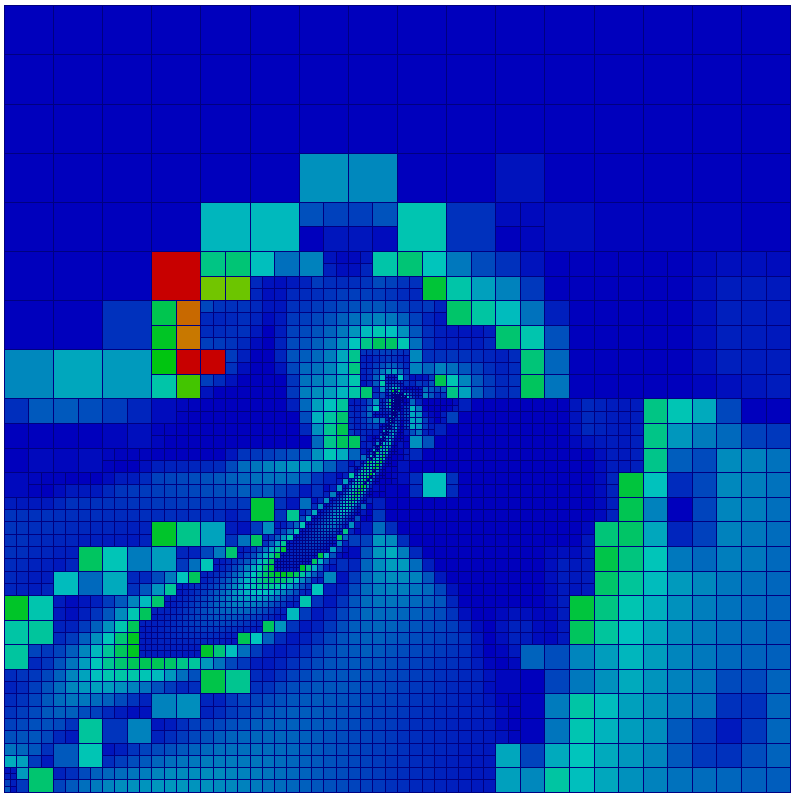}
\end{minipage}
\hfill
\begin{minipage}{0.3\textwidth}
 \includegraphics[width=0.95\textwidth]{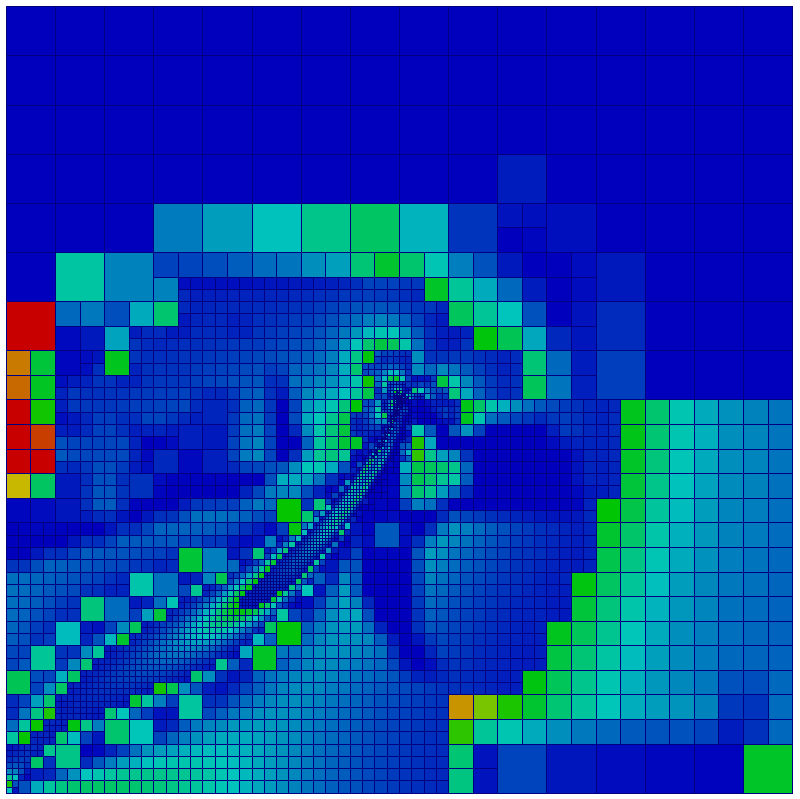}
 \end{minipage}
 \hfill
\begin{minipage}{0.375\textwidth}
 \includegraphics[width=0.95\textwidth]{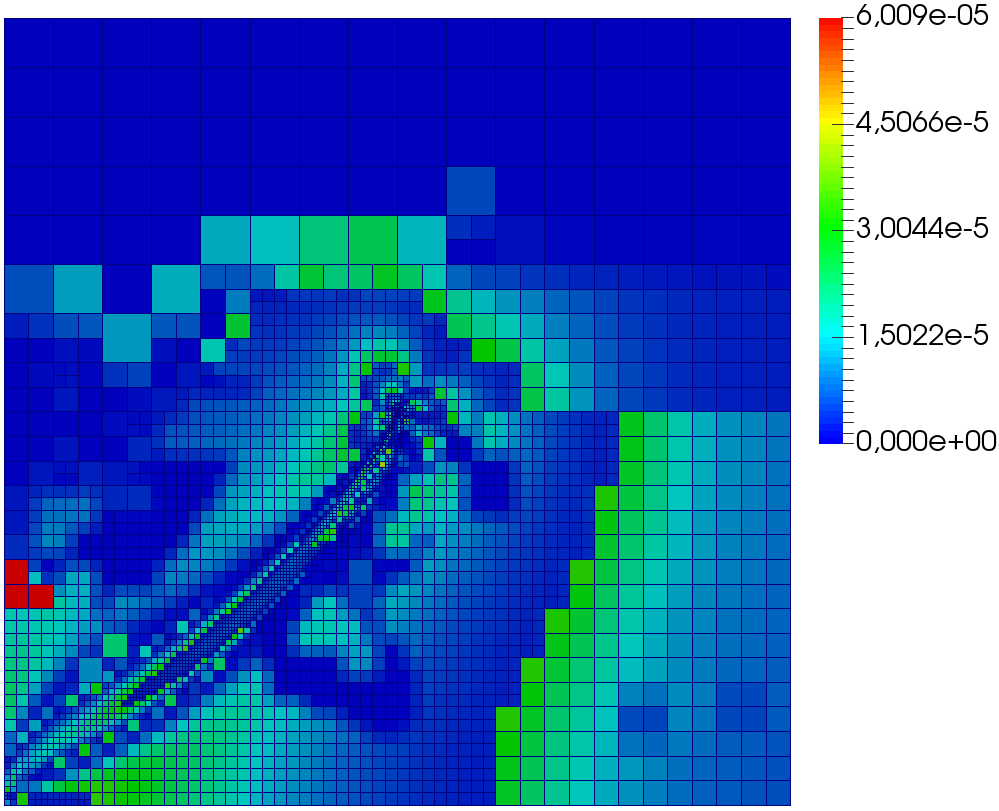}
 \end{minipage}
\caption{The phase-field function and the error indicators, respectively, after six refinement steps depicted in certain time steps (after $0.0116$, $0.0118$ and $0.0125\ \si{s}$) for the single edge notched shear test given on the current adaptive mesh to visualize the refinement strategy.}\label{phasefield_shear_6}
\end{figure}

In Figure~\ref{phasefield_shear_3} and Figure~\ref{phasefield_shear_6}, the phase-field function is depicted after $115$, $119$ and $125$ time steps, respectively. 
The course of the crack to the left lower corner is as expected and as it can
be found in the literature, 
e.g.,~\cite{miehe2010phase,borden2012phase,HeWheWi15,AmGeraLoren15}. 

The snapshots in the first row in Figure~\ref{phasefield_shear_3} and
Figure~\ref{phasefield_shear_6} indicate, that the error estimation and the
corresponding refinement strategy allow to impress the zone around the crack
after three refinement steps; in particular after six steps of adaptive refinement.

\begin{figure}[htbp!]
  \begin{minipage}{0.495\textwidth}
 \includegraphics[width=0.98\textwidth]{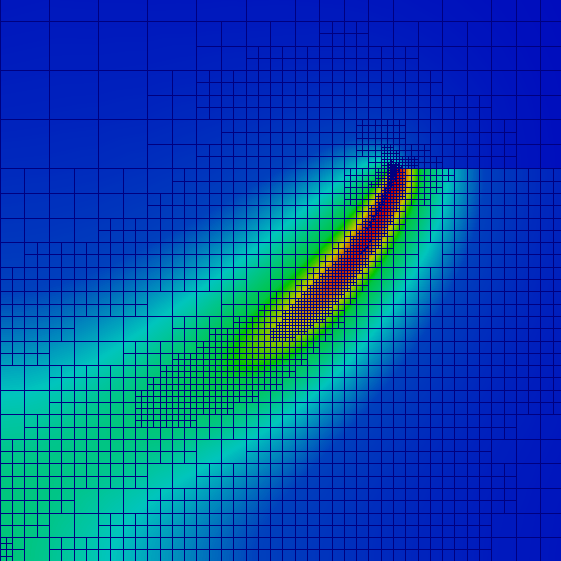}
\end{minipage}
\hfill
\begin{minipage}{0.495\textwidth}
 \includegraphics[width=0.98\textwidth]{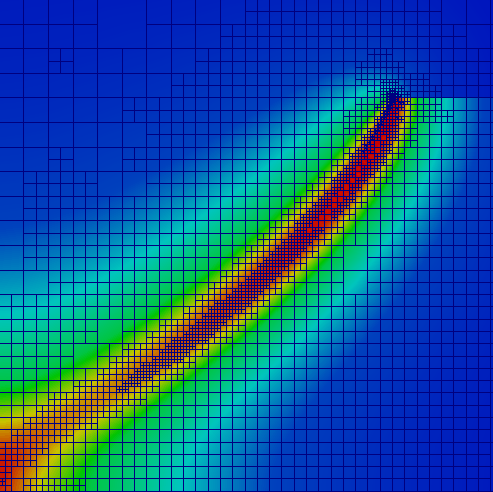}
 \end{minipage}
 \caption{Enhanced extract of the phase-field function after six refinement steps in two certain time steps (after $0.0116$ and $0.0125\ \si{s}$) for the single edge notched shear test given on the current adaptive mesh.}\label{phasefield_shear_6_zoom}
\end{figure}

Further, the snapshots in Figure~\ref{phasefield_shear_6_zoom} allow to see the adaptive mesh especially in the crack domain enlarged in the time steps $116$ and $125$. Mesh cells far away of the fracture are unrefined or very coarse in contrast to the close region of the fracture.

\subsubsection{Results of the single edge notched tension test}

\begin{figure}[htbp!]
\begin{minipage}{0.3\textwidth}
 \includegraphics[width=0.95\textwidth]{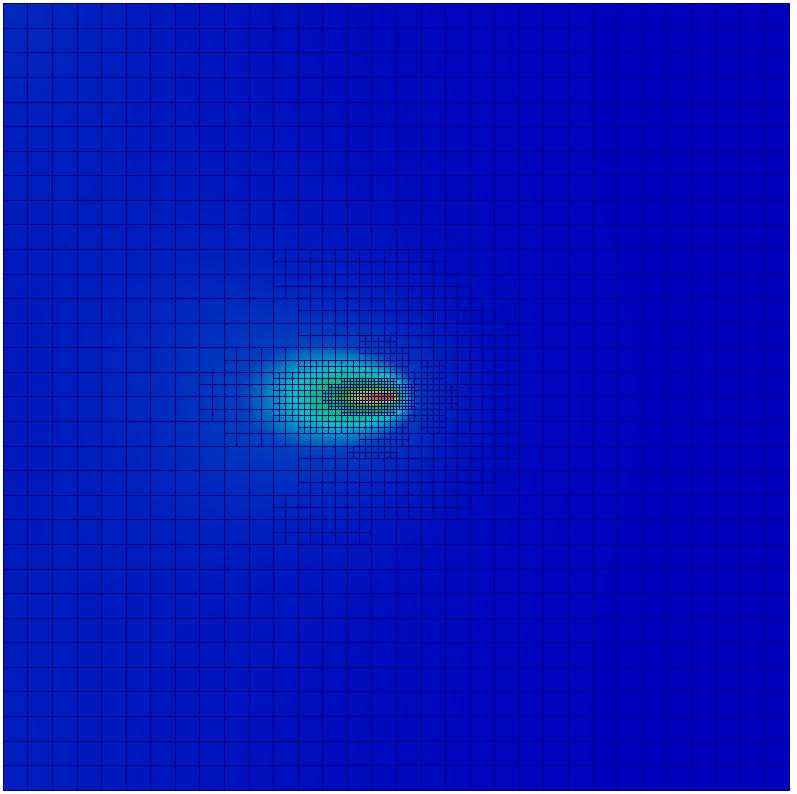}
\end{minipage}
\hfill
\begin{minipage}{0.3\textwidth}
 \includegraphics[width=0.95\textwidth]{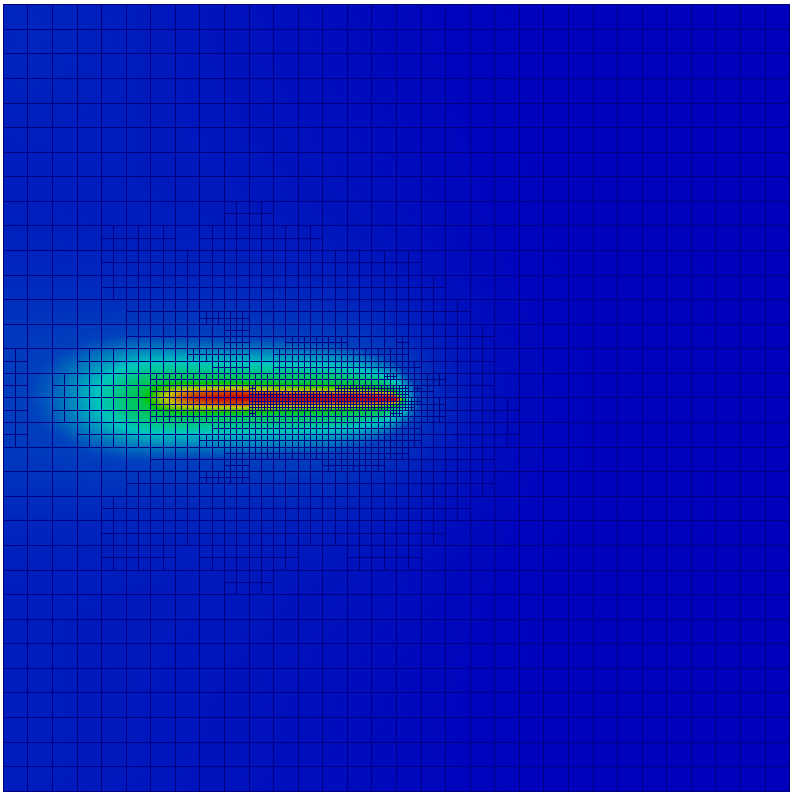}
 \end{minipage}
 \hfill
\begin{minipage}{0.375\textwidth}
 \includegraphics[width=0.95\textwidth]{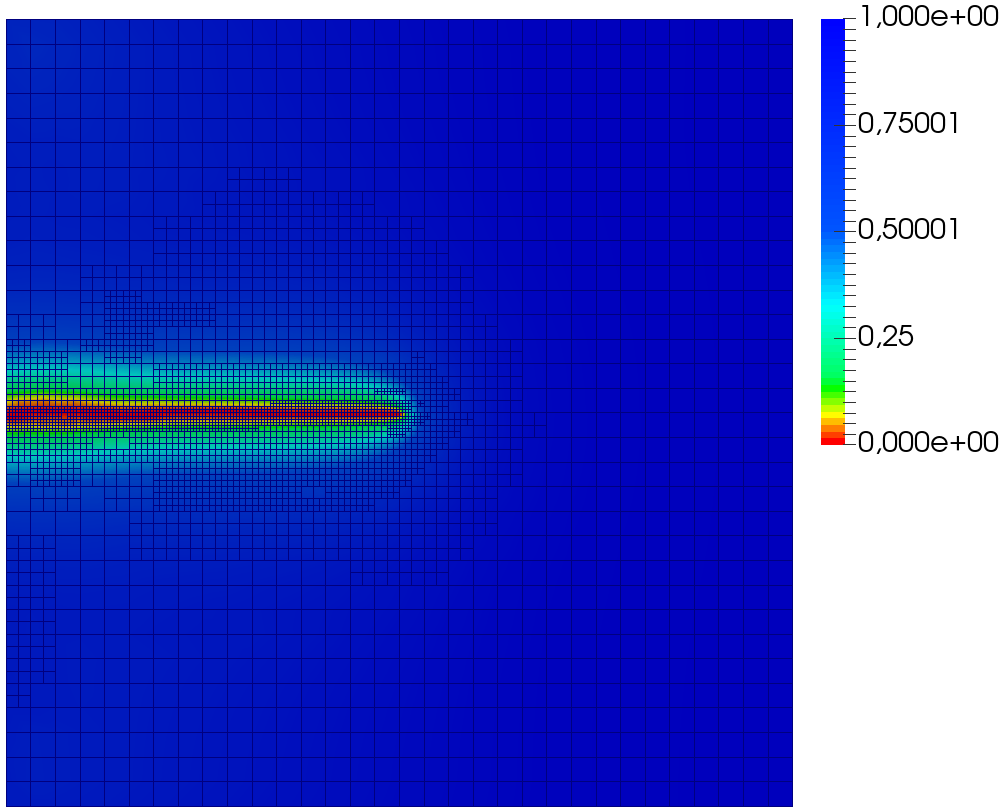}
 \end{minipage}\\
 \vfill
\begin{minipage}{0.3\textwidth}
 \includegraphics[width=0.95\textwidth]{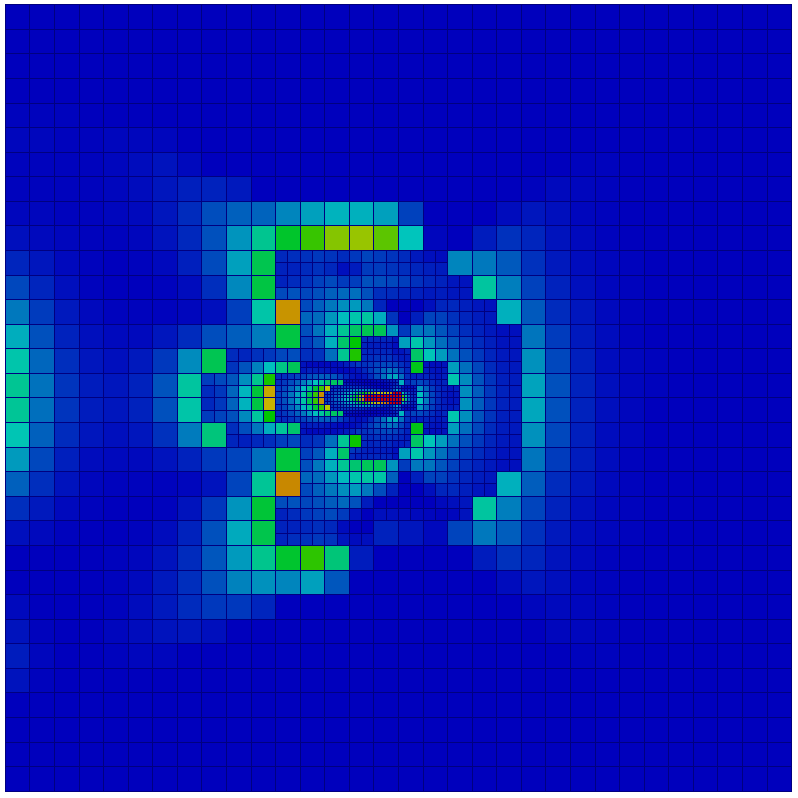}
\end{minipage}
\hfill
\begin{minipage}{0.3\textwidth}
 \includegraphics[width=0.95\textwidth]{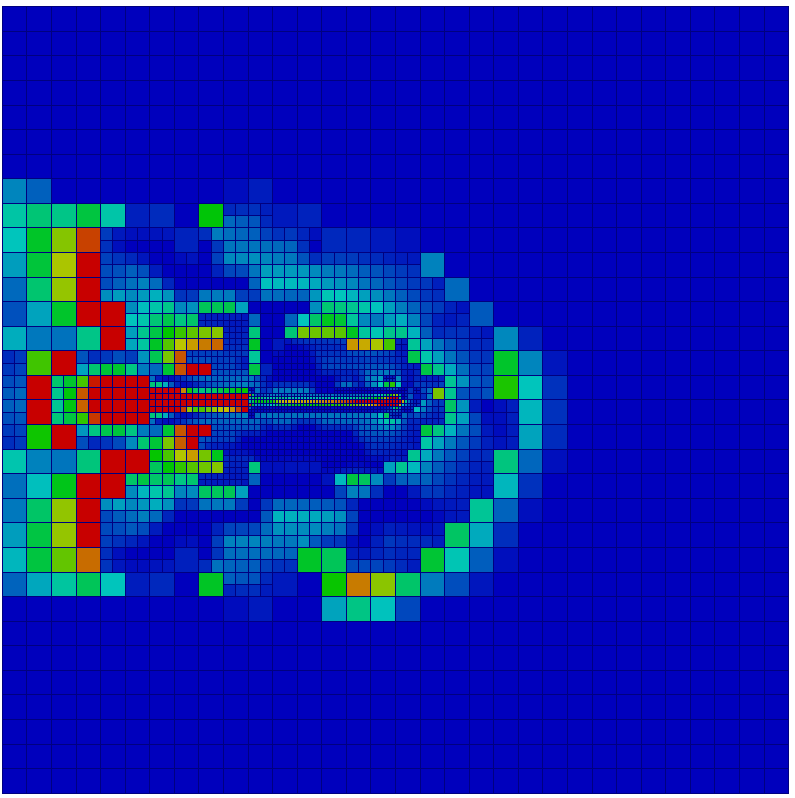}
 \end{minipage}
 \hfill
\begin{minipage}{0.375\textwidth}
 \includegraphics[width=0.95\textwidth]{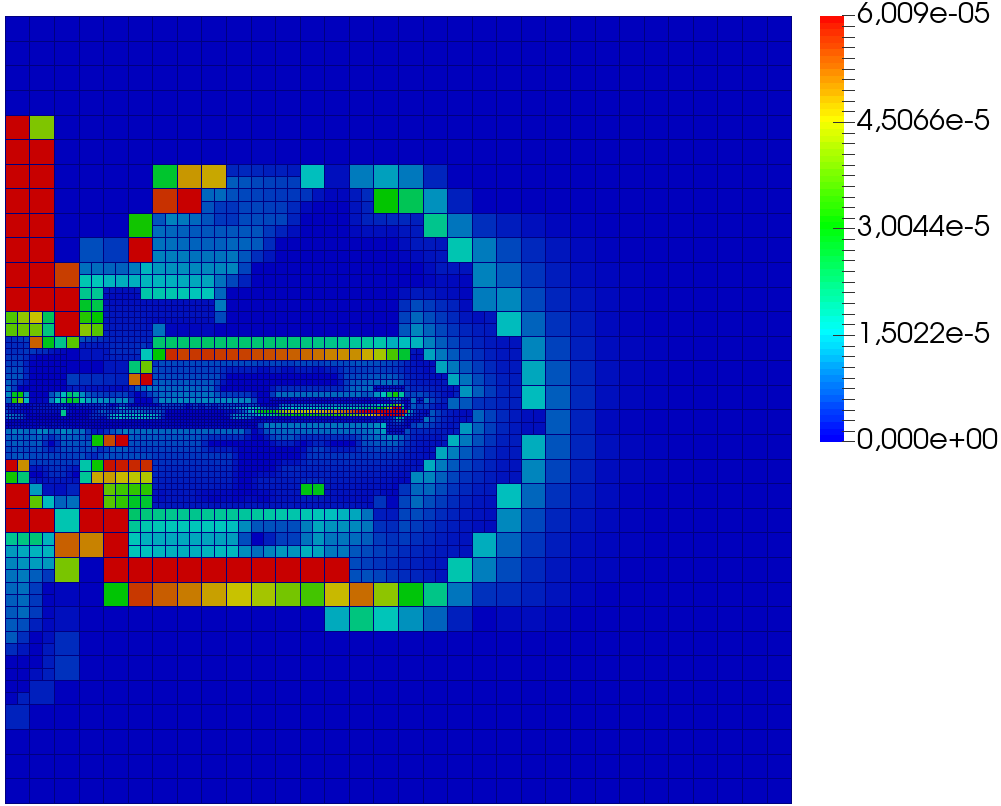}
 \end{minipage}
\caption{The phase-field function and the error indicators, respectively, after three refinement steps in certain time steps (after $0.00615$, $0.00635$ and $0.00676\ \si{s}$) for the single edge notched tension test given on the current adaptive mesh to visualize the refinement strategy.}\label{phasefield_tension_3}
\end{figure}

\begin{figure}[htbp!]
\begin{minipage}{0.3\textwidth}
 \includegraphics[width=0.95\textwidth]{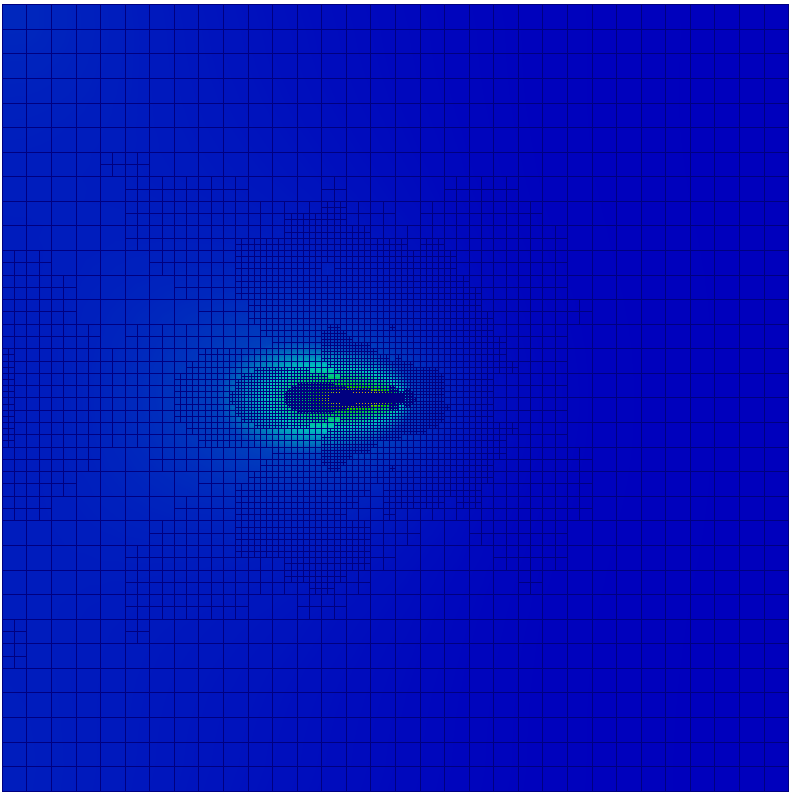}
\end{minipage}
\hfill
\begin{minipage}{0.3\textwidth}
 \includegraphics[width=0.95\textwidth]{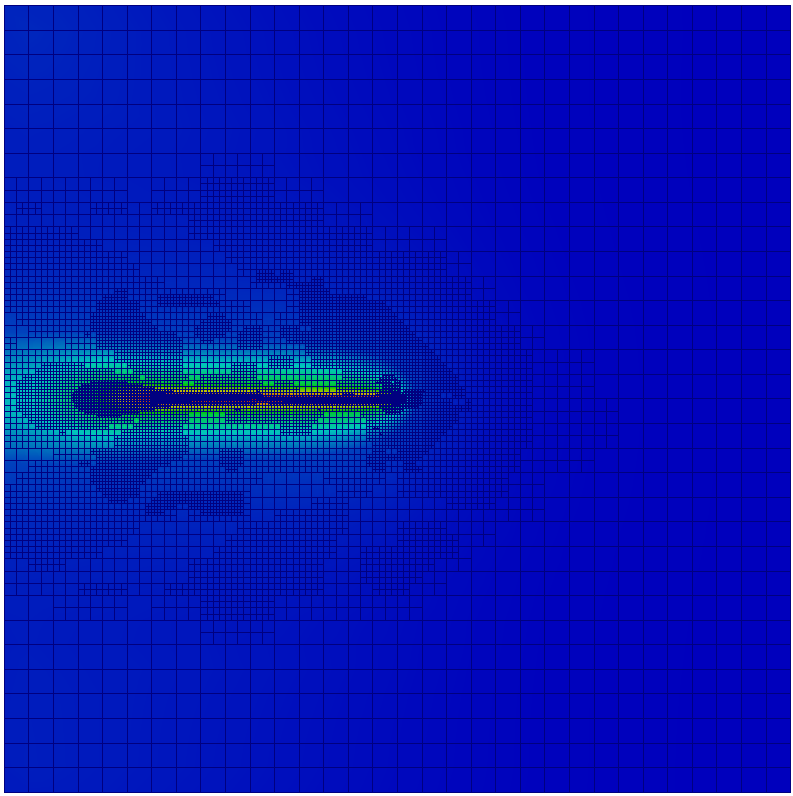}
 \end{minipage}
 \hfill
\begin{minipage}{0.375\textwidth}
 \includegraphics[width=0.95\textwidth]{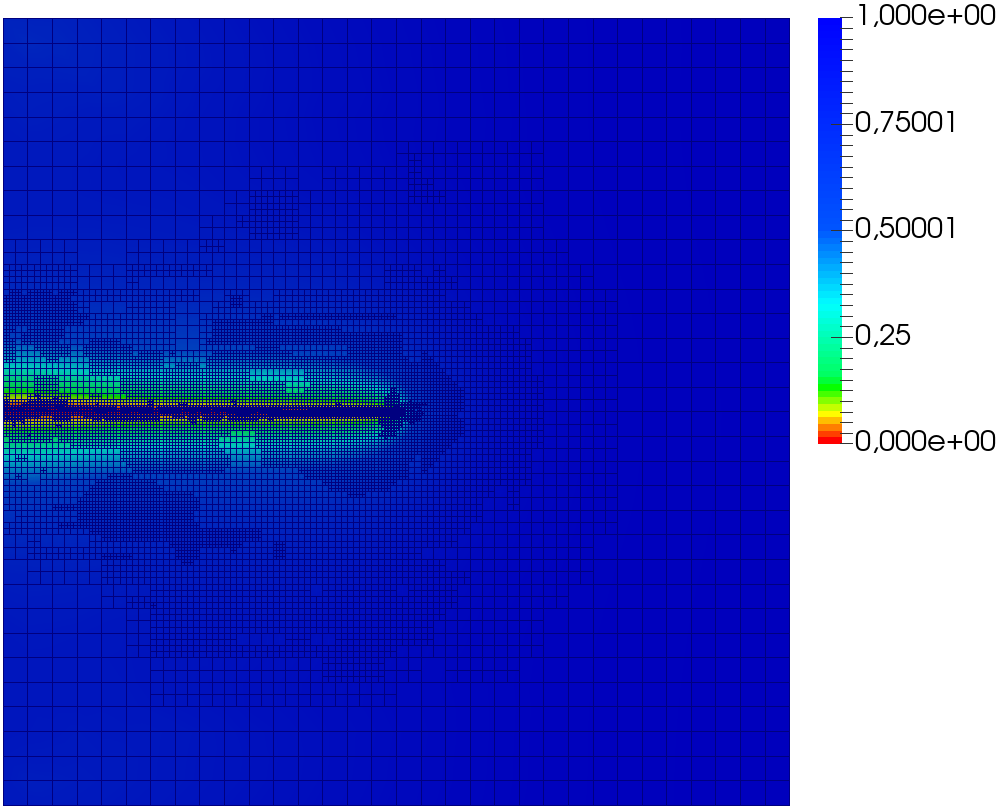}
 \end{minipage}\\
 \vfill
\begin{minipage}{0.3\textwidth}
 \includegraphics[width=0.95\textwidth]{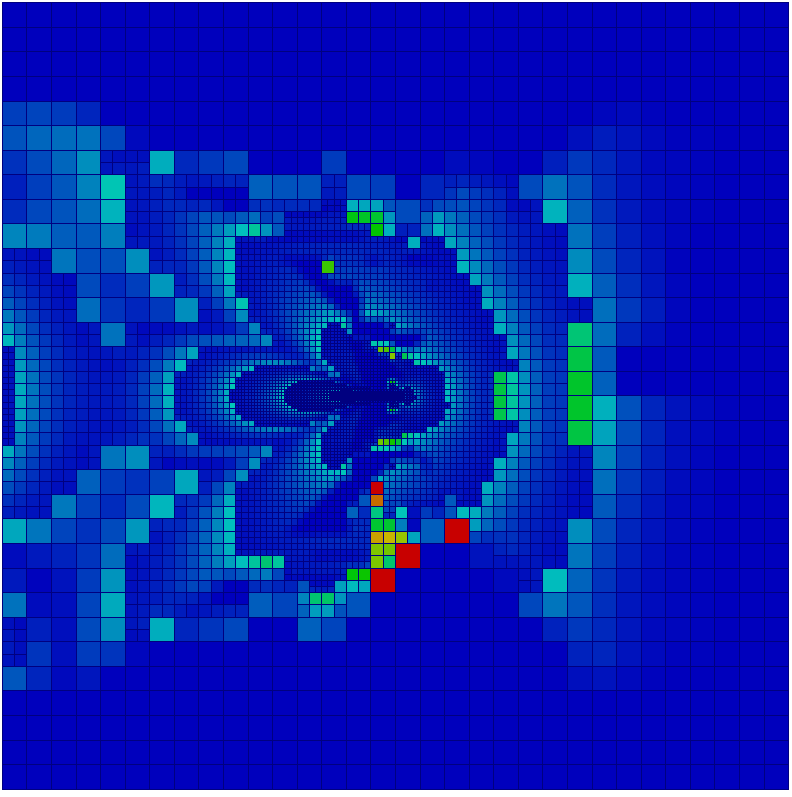}
\end{minipage}
\hfill
\begin{minipage}{0.3\textwidth}
 \includegraphics[width=0.95\textwidth]{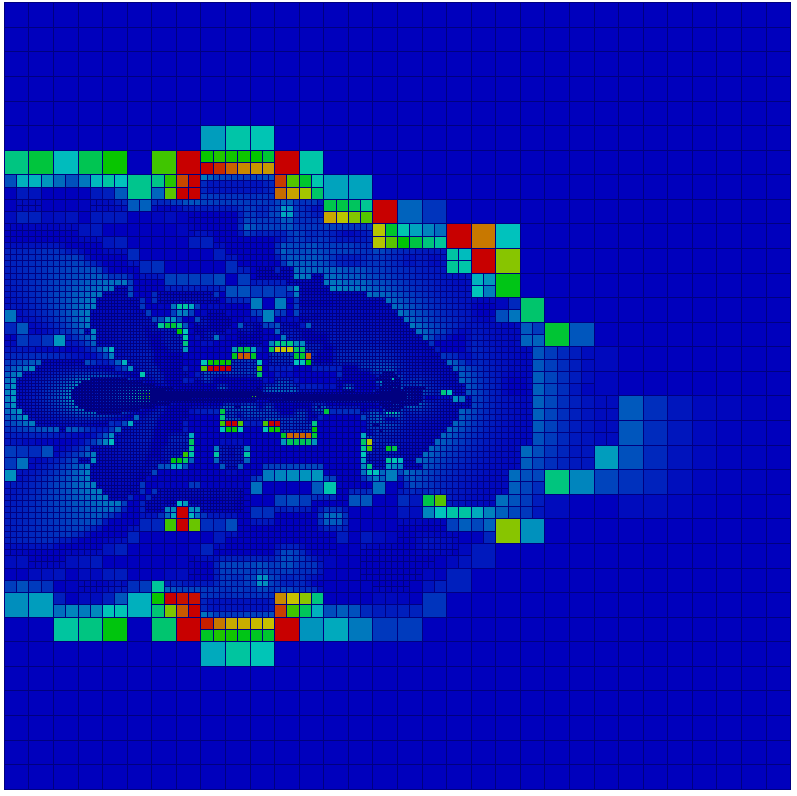}
 \end{minipage}
 \hfill
\begin{minipage}{0.375\textwidth}
 \includegraphics[width=0.95\textwidth]{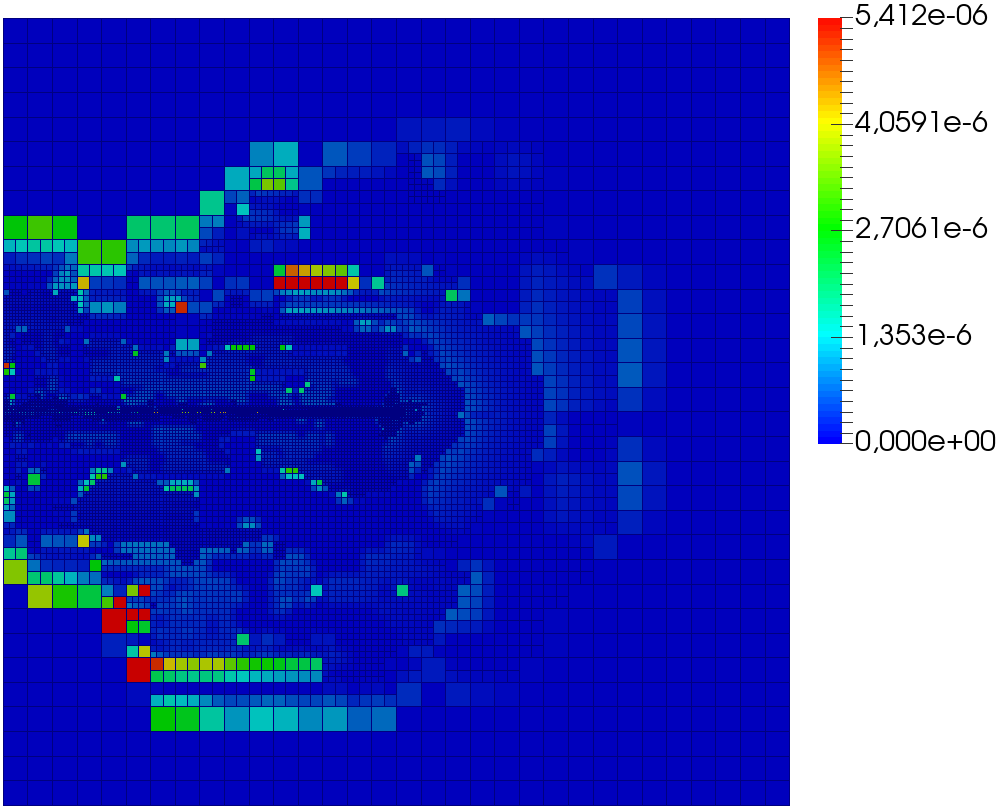}
 \end{minipage}
\caption{The phase-field function and the error indicators, respectively, after six refinement steps in certain time steps (after $0.00615$, $0.00635$ and $0.00674\ \si{s}$) for the single edge notched tension test given on the current adaptive mesh to visualize the refinement strategy.}\label{phasefield_tension_6}
\end{figure}

For the single edge notched tension test, 
in Figure~\ref{phasefield_tension_3} and Figure~\ref{phasefield_tension_6} the phase-field function is displayed after $615$, $635$ and $676$ time steps. 
The course of the crack to the left boundary is as expected and in line with
the literature, 
e.g.,~\cite{miehe2010phase,HeWheWi15}. 

In the Figures~\ref{phasefield_tension_3} and~\ref{phasefield_tension_6}, we observe within the fast crack propagation of the tension test, that the error estimator marks cells with high errors especially in the region before the crack tip, 
which secures, that the crack 
itself moves in a refined region. Further, the plotted error indicators show the symmetry of this test in comparison to the non-symmetric shear test. 

\begin{figure}[htbp!]
  \begin{minipage}{0.495\textwidth}
 \includegraphics[width=0.98\textwidth]{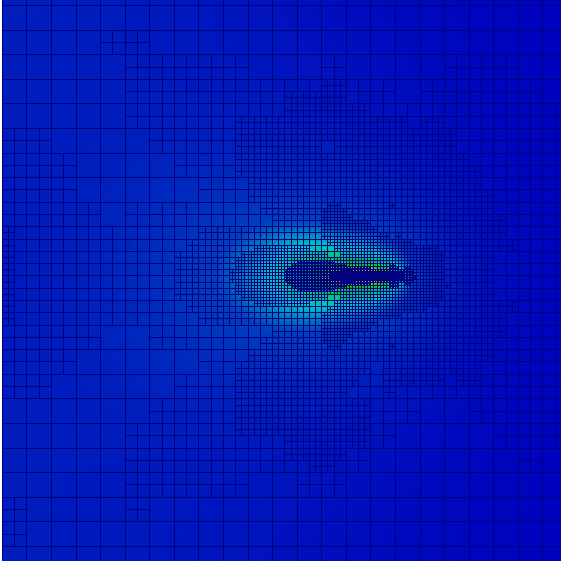}
\end{minipage}
\hfill
\begin{minipage}{0.495\textwidth}
 \includegraphics[width=0.98\textwidth]{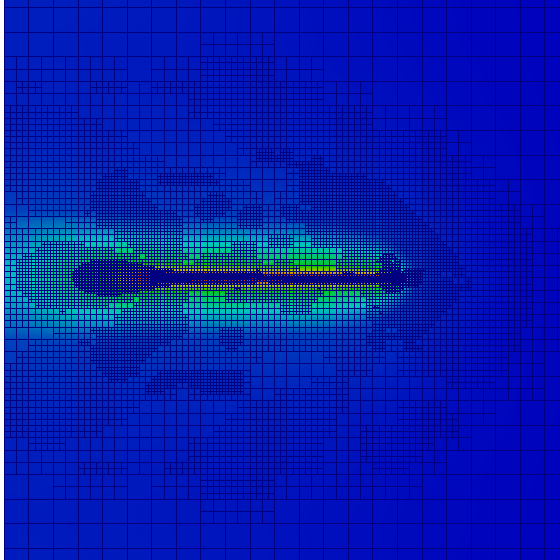}
 \end{minipage}
 \caption{Enhanced extract of the phase-field function after six refinement steps in certain time steps (after $0.00615$ and $0.00635\ \si{s}$) for the single edge notched tension test given on the current adaptive mesh.}\label{phasefield_tension_6_zoom}
\end{figure}

To improve the transparency of the new developed error estimator and the adaptive refinement strategy regarding the detected error, 
Figure~\ref{phasefield_tension_6_zoom} contains zoomed snapshots of the phase field function and the current adaptive meshes are given at the time steps $615$ and $635$. 


\section{Conclusions}
\label{sec_conclusions}

The focus of this work was to develop a residual-type error estimator for
phase-field fracture propagation problems. Due to the fracture irreversibility
constraint, we deal with a variational inequality in time. 
Based on these theoretical advancements, we developed an adaptive solution
strategy for the monolithically-coupled displacement/phase-field system.
We investigated the performance by the help of two numerical
configurations. First, we considered the so-called single edge notched 
shear test in which a curved fracture develops. Fixing the phase-field 
regularization parameter and varying the spatial mesh parameter, we obtained 
excellent convergence behavior of the load-displacement curves. The same 
observations were made for the evolution of the bulk and the crack energy. In view of mesh
adaptivity, we obtained localized mesh refinement in the (a priori unknown)
fracture region. 
For the second numerical example, the single edge notched
tension test, we noticed that here we have very fast, brutal crack
growth, which is challenging for mesh refinement strategies. Again, we 
observed very convincing findings. In ongoing work, we will apply the proposed 
residual-based error estimator to a phase-field fracture model in incompressible
solids. Furthermore, we will provide the proofs of  
the reliability and efficiency of the proposed estimators, with constants independent of
the chosen parameter $\epsilon$.


\section*{Acknowledgments}
This work has been supported by the German Research Foundation, Priority Program 1748 (DFG SPP 1748) named
\textit{Reliable Simulation Techniques in Solid Mechanics. Development of
Non-standard Discretization Methods, Mechanical and Mathematical
Analysis}. Our subproject within the SPP1748 reads \textit{Structure
  Preserving Adaptive Enriched Galerkin Methods for Pressure-Driven 3D
  Fracture Phase-Field Models} 
(WA 4200/1-1 and WI 4367/2-1 and WO 1936/5-1).


\bibliographystyle{abbrv}
\bibliography{lit}

\end{document}